\documentclass[a4paper,12pt]{article}

\usepackage{amsmath}
\usepackage{amssymb}
\usepackage{url}

\ifx\pdfoutput\undefined
\usepackage{graphicx}
\usepackage{epsfig}

\else
\usepackage[pdftex]{graphicx}

\fi

\usepackage{color}

\newtheorem{thm}{Theorem}[subsection]
\newtheorem{lem}[thm]{Lemma}
\newtheorem{prp}[thm]{Proposition}
\newtheorem{cor}[thm]{Corollary}

\newtheorem{prob}[thm]{Problem}
\newenvironment{proof}{{\it Proof}.\ }{\hfill$\square$\par\medskip}

\def\cF{{\mathcal F}}

\def\cH{{\mathcal H}}

\def\cT{{\mathcal T}}

\def\frD{{\mathfrak D}}

\def\frM{{\mathfrak M}}

\def\frd{{\mathfrak d}}

\def\frh{{\mathfrak h}}

\def\frm{{\mathfrak m}}

\newcommand\RR{{\mathbb R}}
\newcommand\Sph{{\mathbb S}}
\newcommand\ZZ{{\mathbb Z}}

\newcommand\SetOf[2]{\bigl\{#1\,\bigm|\,#2\bigr\}}
\newcommand\id{\operatorname{id}}

\def\Sym{\operatorname{Sym}}

\def\SetOf#1#2{\left\{\left.#1\vphantom{#2}\ \right|\ #2\vphantom{#1}\right\}}

\def\to{\rightarrow}

\def\conv{\operatorname{conv}}
\def\cone{\operatorname{cone}}
\def\core{\operatorname{core}}
\def\mod{\operatorname{mod}}

\def\phi{\varphi}
\def\epsilon{\varepsilon}
\def\int{\operatorname{int}}
\def\codim{\operatorname{codim}}
\def\ost{\operatorname{stint}}
\def\st{\operatorname{st}}
\def\lk{\operatorname{lk}}

\def\relint{\operatorname{relint}}
\def\ind{\operatorname{ind}}

\newcommand\DD{{\mathbb D}}
\def\nothing{\rule{0mm}{1ex}}

\def\oK{{\overline{K}}}

\def\unfolding#1{{\widetilde{#1}}}
\def\partunfold#1{{\widehat{#1}}}
\def\tK{\unfolding{K}}
\def\odd#1{{#1}_{\mathrm{odd}}}

\title{Branched Coverings, Triangulations,\\ and $3$-Manifolds}
\author{Ivan Izmestiev and Michael Joswig}
\date{\today}

\begin{document}
\maketitle

\begin{abstract}
  A canonical branched covering over each sufficiently good simplicial complex is constructed.  Its structure depends on
  the combinatorial type of the complex.  In this way, each closed orientable $3$-manifold arises as a branched covering
  over~$\Sph^3$ from some triangulation of~$\Sph^3$.  This result is related to a theorem of Hilden~\cite{Hil} and
  Montesinos~\cite{Mon}.  The branched coverings introduced admit a rich theory in which the group of projectivities,
  defined in~\cite{ProjSimplePolytopes}, plays a central role.
\end{abstract}

\section{Introduction}\label{sec:intro}

A celebrated theorem of Hilden and Montesinos says that each oriented $3$-manifold can be obtained as a special branched
covering space over the $3$-sphere; for a precise formulation see Theorem~\ref{thm:Hilden-Montesinos} below.  The
purpose of this paper is to show how these (and many other) branched coverings can be described in a purely
combinatorial way.

There is already quite an extensive literature on the combinatorial treatment of branched coverings of manifolds.  Often
this work is restricted to surfaces; for instance see Gross and Tucker~\cite{GT}.  This is analogous to the historical
development of the topological theory of branched coverings.  It has its roots in the theory of Riemann surfaces and
later extended to higher-dimensional manifolds.  However, the foundations of the theory of branched coverings for a
wider class of topological spaces were laid only in the 1950s by Fox~\cite{Fox}.  One of the few attempts to give a
combinatorial treatment of a general class of branched coverings is due to Mohar~\cite{Moh}.  He applied voltage graphs,
see~\cite{GT}, and Fox's theory to obtain an encoding of branched coverings of (pseudo-)simplicial complexes.

Our point of view is a different one.  We focus on the explicit construction of a very special class of branched
coverings, which we call \emph{unfoldings}.  As a key property these unfoldings are canonically associated to a
triangulation of the base space.  Phrased differently, starting from a (sufficiently good) triangulation of a
topological space we give an elementary combinatorial description of the branched covering space of an unfolding.  This
is remarkable since surprisingly many branched covering maps can be described in this way.  In particular, all branched
covering maps in the aforementioned theorem of Hilden and Montesinos arise.

The key tool for our investigation is the group of projectivities~$\Pi(K)$ of a finite simplicial complex~$K$, which has
been explored in~\cite{ProjSimplePolytopes}.  Originally devised for the study of certain coloring problems this group
turns out to behave similar to a fundamental group, whereas the \emph{complete unfolding} plays the role of the
universal covering.  In particular, the group of projectivities~$\Pi(\unfolding{K})$ of the complete unfolding is always
trivial.

It is essential that the unfoldings depend on the combinatorial properties of~$K$.  Although an arbitrary subdivision
of~$K$ does not change the PL-type, it can influence the group of projectivities and the unfoldings in a rather
unpredictable way.  On the other hand, in order to prove Theorem~\ref{thm:main-result} we make use of a variant of the
Simplicial Approximation Theorem.  This requires a special type of subdivision which preserves the group of
projectivities and yields an equivalent unfolding.  We give an explicit construction of such a subdivision which we call
the \emph{anti-prismatic subdivision}.

It should be pointed out that our results could also be stated in the language of voltage graphs.  Since our proofs,
however, seem to require very different techniques we leave this to the interested reader.

We outline the organization of the paper.

We start by recalling the definition and the basic properties of the group of projectivities.  Then we construct the
complete and partial unfoldings of an arbitrary pure simplicial complex.  Here a technical difficulty arises: In
general, an unfolding may have a more complicated structure than a simplicial complex.  In the literature objects of
this class are often called pseudo-simplicial complexes.  However, we show that this is only a minor problem.  Firstly,
one can extend the notion of a projectivity to pseudo-simplicial complexes.  Secondly, after an anti-prismatic
subdivision the complete unfolding becomes a simplicial complex.  The technical details are deferred to the Appendix.

The next section is devoted to a more thorough investigation of the unfoldings.  From the theory of coverings it is
familiar that certain local connectivity properties are required in order to yield a satisfying theory.  In a similar
way, we introduce additional restrictions on the local structure of the complex.  The crucial property of these
\emph{nice} complexes is that their dual block structure is good enough.  The class of nice complexes includes all
PL-manifolds as well as all (locally finite) graphs.  It turns out that one can find a system of generators for the
group of projectivities of a nice complex.  This directly generalizes the corresponding
result~\cite[Theorem~8]{ProjSimplePolytopes} on PL-manifolds.

In Section~\ref{sec:BrCov} we briefly recall Fox's theory of branched coverings~\cite{Fox}.  Then we prove
Theorem~\ref{thm:UnfoldingIsBranched}: The unfoldings of nice complexes are, in fact, branched coverings.  The branch
set of the complete unfolding is the \emph{odd subcomplex}, formed by the codimension-$2$-faces whose links are
non-bipartite graphs.  Moreover, the complete unfolding~$\unfolding{K}$ is regular, and the group of
projectivities~$\Pi(K)$ is its group of covering transformations.  Besides, we show that the complete unfolding is the
regularization of the partial unfolding.

The final section contains a discussion of the unfoldings of PL-manifolds.  In particular, we study the problem to
determine which branched coverings of a PL-manifold arise as unfoldings.  The proof of our key
result~\ref{thm:main-result} can be sketched as follows.  For a given closed oriented $3$-manifold~$M$ we start with a
branched covering $f:M\to\Sph^3$ as in the Hilden-Montesinos Theorem.  The covering map~$f$ is branched over a
knot~$L\subset\Sph^3$.  Up to equivalence it is characterized by its monodromy homomorphism~$\frm_f:\pi(\Sph^3\setminus
L)\to S_3$.  Then we construct a triangulation of the pair $(\Sph^3,L)$, where the group of projectivities is isomorphic
to~$S_3$, the symmetric group of degree~$3$, and $L$ is the odd subcomplex.  It follows that $M$ is the partial
unfolding of this triangulation.

\section{Projectivities and the unfolding of a\\ simplicial complex} 

\subsection{The group of projectivities}

Throughout the whole paper let $K$ be a $d$-dimensional locally finite simplicial complex.  Moreover, we assume that $K$
is \emph{pure}, that is, each face of $K$ is contained in a face of dimension $d$.  The $d$-dimensional faces of $K$ are
called \emph{facets}, the faces of codimension~$1$ are called \emph{ridges}.  The dual graph~$\Gamma(K)$ has the facets
of~$K$ as nodes, and an edge connects two such nodes if the corresponding facets share a common ridge.  Suppose that a
ridge~$\rho$ is contained in two facets $\sigma$ and~$\tau$.  Then there is unique vertex~$v(\sigma,\tau)$ of~$\sigma$,
which is not contained in~$\tau$.  We denote the set of vertices of~$\sigma$ by~$V(\sigma)$, and we introduce a
bijective map

$$\langle\sigma,\tau\rangle:V(\sigma)\to V(\tau):w\mapsto\left\{\begin{array}{cl} v(\tau,\sigma) & \text{if
      $w=v(\sigma,\tau)$,}\\ w & \text{otherwise.}\end{array}\right.$$

The map $\langle\sigma,\tau\rangle$ is called the \emph{perspectivity} from $\sigma$ to~$\tau$.  A \emph{facet path}
in~$K$ is a sequence $\gamma = (\sigma_0,\sigma_1,\ldots,\sigma_n)$ such that $\sigma_0,\sigma_1,\ldots,\sigma_n$ are
facets and two consecutive facets $\sigma_i$ and~$\sigma_{i+1}$ are neighbors in~$\Gamma(K)$ for all $0\le i<n$.  Now
the \emph{projectivity} along~$\gamma$ is defined as the product of perspectivities

$$\langle\gamma\rangle = \langle\sigma_0,\sigma_1,\ldots,\sigma_n\rangle = \langle\sigma_0,\sigma_1\rangle \cdots
\langle\sigma_{n-1},\sigma_n\rangle.$$

\begin{figure}[htbp]
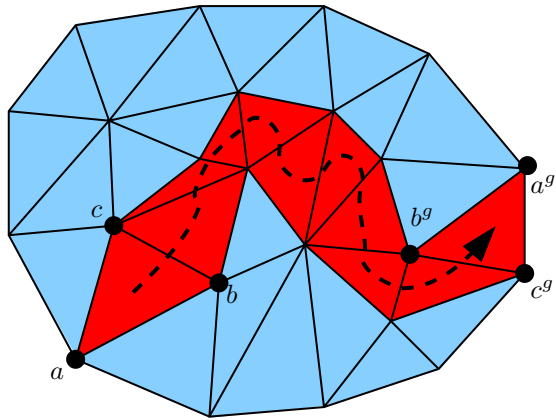

  \begin{center}
    \input path-proj.pstex_t 
    \caption{Facet path~$\gamma$ and projectivity $g=\langle\gamma\rangle$.}
    \label{fig:path-proj}
  \end{center}
\end{figure}

The \emph{inverse path} of~$\gamma$ is denoted by $\gamma^- = (\sigma_n,\sigma_{n-1},\ldots,\sigma_0)$.  We write
$\gamma\eta$ for the concatenation of~$\gamma$ with some facet path~$\eta=(\sigma_n,\ldots,\sigma_m)$.  The facet path
$\gamma$ is \emph{closed} if $\sigma_0=\sigma_n$.  In this case we also call $\gamma$ a \emph{facet loop} based
at~$\sigma_0$.  The set of projectivities along facet loops based at~$\sigma_0$ forms a group, the \emph{group of
  projectivities} at~$\sigma_0$, which is written as $\Pi(K,\sigma_0)$.  The group $\Pi(K,\sigma_0)$ is a subgroup of
the group $\Sym(V(\sigma_0))$ of all permutations of~$V(\sigma_0)$.  If $\sigma$ and~$\tau$ are facets of~$K$ which can
be joined by a facet path, that is, they are contained in the same connected component of~$\Gamma(K)$, then
$\Pi(K,\sigma)$ is isomorphic to~$\Pi(K,\tau)$ as a permutation group; equivalently, the groups become conjugate after
an arbitrary identification between the sets $V(\sigma)$ and $V(\tau)$.  In particular, if $K$ is \emph{strongly
  connected}, that is, the dual graph~$\Gamma(K)$ is connected, this yields a subgroup~$\Pi(K)$ of the symmetric group
$S_{d+1}$ of degree~$d+1$, which is well defined up to conjugation.

Groups of projectivities of simplicial complexes have been introduced in~\cite{ProjSimplePolytopes}.

We will use two alternative notations $fg$ and $g\circ f$ for the composition of maps $f:X\to Y$ and $g:Y\to Z$.  The
first notation is used in the context of projectivities, while the second is used in all other cases.  The
projectivities operate on the right. Throughout we use the same notation for a simplicial complex and its geometric
realization.

Occasionally, we want to examine topological properties of facet paths.  Observe that each facet
path~$\gamma=(\sigma_0,\ldots,\sigma_n)$ in~$K$ induces a piecewise linear path $\overline\gamma$ in the geometric
realization of~$K$: Join the barycenter of each facet~$\sigma_i$ by linear paths to the barycenters of the common ridges
$\sigma_i\cap\sigma_{i-1}$ and $\sigma_i\cap\sigma_{i+1}$ of the neighboring facets $\sigma_{i-1}$ and~$\sigma_{i+1}$,
respectively.  The facet path~$\gamma$ is closed if and only if the induced piecewise linear path~$\overline\gamma$ is
closed.  Often we identify $\gamma$ with~$\overline\gamma$.  Moreover, we write $[\gamma]$ for the homotopy class
of~$\overline\gamma$ with endpoints fixed.

By $\Pi_0(K,\sigma_0)$ denote the subgroup of~$\Pi(K,\sigma_0)$ of projectivities along facet loops which are
null-homotopic.  We call $\Pi_0(K,\sigma_0)$ the \emph{reduced group of projectivities}.

\begin{prp}\label{prp:reduced-normal}
  The group~$\Pi_0(K,\sigma_0)$ is a normal subgroup of~$\Pi(K,\sigma_0)$.
\end{prp}

\begin{proof}
  Let $\gamma$ and~$\eta$ be facet loops based at~$\sigma_0$, and suppose that $\eta$ is null-homotopic.  Then the facet
  loop $\gamma^-\eta\gamma$ is also null-homotopic.
\end{proof}

A simplicial map is \emph{non-degenerate} if it takes each simplex to a simplex of the same dimension.

\begin{prp}\label{prp:non-degenerate}
  Let $f:K\to L$ be a non-degenerate simplicial map between pure complexes of the same dimension.  Then for each pair of
  facets $\sigma_0\in K$ and $\tau_0\in L$ such that $f(\sigma_0)=\tau_0$ there is a canonical homomorphism
  $f_*:\Pi(K,\sigma_0)\to\Pi(L,\tau_0)$.  Moreover, $f_*$ is injective.
\end{prp}

\begin{proof}
  Note that the images of neighboring facets under $f$ either coincide or are neighboring facets as well. Hence, for any
  facet path $\gamma$ in $K$, we can form a facet path $f_\sharp(\gamma)$ in $L$ deleting from the sequence of images of
  facets in $\gamma$ every term that coincides with the preceding one. Clearly, $\langle f_\sharp(\gamma)\rangle =
  f^{-1}\langle\gamma\rangle f$, where $f^{-1}$ is considered as a map from $V(\tau_0)$ to~$V(\sigma_0)$.  In
  particular, the projectivity $\langle f_\sharp(\gamma)\rangle$ is the identity if and only if $\langle\gamma\rangle$
  is the identity.  This implies that the map $f_*:\Pi(K,\sigma_0)\to\Pi(L,\tau_0),\, \langle\gamma\rangle\mapsto\langle
  f_\sharp(\gamma)\rangle$ is a well defined monomorphism.
\end{proof}

Thus the group of projectivities provides an obstruction for the existence of a non-degenerate simplicial map between
simplicial complexes.  Letting $L$ be the $d$-simplex where $d=\dim K$ this result implies that if the vertices of~$K$
can be properly colored with $d+1$~colors, then $\Pi(K)=1$, see~\cite[Proposition~6]{ProjSimplePolytopes}.

\subsection{The unfoldings}\label{subsec:unfolding}

Here we introduce two geometric objects defined by the combinatorial structure of the simplicial complex $K$. These are
the complete unfolding~$\unfolding{K}$ and the partial unfolding~$\partunfold{K}$ together with canonical maps
$p:\unfolding{K} \to K$ and $r:\partunfold{K} \to K$. The spaces $\unfolding{K}$ and $\partunfold K$ arise as special
quotient spaces of collections of geometric simplices.  However, they may not be simplicial complexes.  Therefore, we
first have to introduce a slightly more general concept.

Let $\Sigma$ be a collection of pairwise disjoint geometric simplices.  We assume that we are given attaching data of
the following form.  For some pairs of simplices~$\sigma$ and~$\tau$ we have a simplicial isomorphism from a subcomplex
of~$\sigma$ to a subcomplex of~$\tau$.  By performing the corresponding identifications in an arbitrary order we obtain
a quotient space $\Sigma/\!\sim$.  Suppose that for each simplex $\sigma\in\Sigma$ the restriction of the quotient map
$\Sigma\to\Sigma/\!\sim$ to $\sigma$ is bijective (that is, within each simplex there are no self-identifications).
Following Hilton and Wylie~\cite{163.17803} we call $\Sigma/\!\sim$ a \emph{pseudo-simplicial complex} or, shortly, a
\emph{pseudo-complex}.  Observe that in a pseudo-complex the intersection of two simplices is not necessarily a single
simplex.  A map between two pseudo-simplicial complexes is called \emph{simplicial} if it takes each simplex of the
first linearly to a simplex of the second.  As an example of a pseudo-simplicial complex consider two copies of the
$d$-dimensional simplex identified along the boundary.  Clearly, the result is homeomorphic to the $d$-sphere~$\Sph^d$.
It is easily seen that the barycentric subdivision of a pseudo-simplicial complex is a simplicial complex.  In
particular, a pseudo-complex has a natural PL-structure.

Similar to a simplicial complex, a pseudo-complex also has a dual graph, which may have multiple edges between nodes but
no loops.  The concepts of facet paths, perspectivities, and projectivities carry over.  The only difference is that, in
a facet path, it is necessary to specify the ridges between the facets.  Other than that everything discussed so far
also holds for projectivities in pseudo-complexes.  We omit the details.  Besides, in the
Appendix~\ref{subsec:anti-prismatic} we construct a subdivision of a pseudo-complex which is a simplicial complex, and
which does not change the group of projectivities.

For the rest of the section let $K$ be strongly connected with a fixed facet~$\sigma_0$.

Consider the disjoint union~$\Sigma(K)$ of facets of~$K$ and the product~$\oK=\Sigma(K)\times\Pi(K,\sigma_0)$.  Each
pair~$(\sigma,g)$ is a copy of the geometric $d$-simplex~$\sigma$.  Thus we have a set of natural affine
isomorphisms~$(\sigma,g)\to\sigma$ which induce the projection $\oK\to K$.  We glue the simplices~$(\sigma,g)$ as
follows.  For each facet $\sigma$ of~$K$ choose some facet path $\gamma_\sigma$ from $\sigma_0$ to~$\sigma$.  Suppose
that $\rho$ is a common ridge in~$K$ of the facets~$\sigma$ and $\tau$.  Then we glue $(\sigma,g)$ and $(\tau,h)$ with
respect to the affine map induced by the identity map on~$\rho$ if the equation
\begin{equation}\label{eq:gluing}
gh^{-1} = \langle\gamma_\sigma\rangle \langle\sigma,\tau\rangle \langle{\gamma_\tau}^-\rangle
\end{equation}
holds in~$\Pi$.  Let $\sim$ be the equivalence relation generated by this gluing
strategy.  The resulting pseudo-simplicial complex
$$\tK=\oK/\!\!\sim$$
is called the \emph{complete unfolding} of~$K$. The \emph{complete unfolding map}
$p:\unfolding{K}\to K$ factors the projection $\oK\to K$ in a natural way.

In order to facilitate the investigation we give an alternative description of the complete unfolding.  Fix a
coloring~$b_0:V(\sigma_0)\to\{0,\ldots,d\}$ of the vertices of the base facet~$\sigma_0$.  For any facet path~$\eta$
from~$\sigma_0$ to some facet~$\sigma$ we obtain an induced coloring $\langle\eta^-\rangle b_0$ of~$V(\sigma)$.  We call
such a coloring of~$V(\sigma)$ \emph{admissible}.  The admissible colorings of the vertex set of a fixed facet
correspond to the elements of the group of projectivities.  We consider the disjoint union~$\oK'$ of
simplices~$(\sigma,b)$, where $b$ is an admissible coloring of the facet~$\sigma$.  Now we glue $(\sigma,b)$ and
$(\tau,c)$ with respect to the identity map on the common ridge $\rho$ of~$\sigma$ and~$\tau$ provided that the
respective restrictions of the colorings $b$ and $c$ to the ridge~$\rho$ coincide.  As the quotient we obtain a
pseudo-simplicial complex~$\unfolding{K}'$.

\begin{prp}
  The two constructions above are simplicially equivalent, that is, there exists a simplicial isomorphism between
  $\unfolding{K}$ and~$\unfolding{K}'$ which commutes with the canonical projections to~$K$.  Moreover, the
  combinatorial structure of the complete unfolding~$\unfolding{K}$ neither depends on the choice of the
  facet~$\sigma_0$ nor on the choice of the facet paths~$\gamma_\sigma$ nor on the choice of the coloring~$b_0$.
\end{prp}

\begin{proof}
  The collection~$\oK'$ of admissibly colored facets is isomorphic to~$\oK$ by virtue of the map
  $$\iota:(\sigma,\langle\eta^-\rangle b_0)\mapsto(\sigma,\langle\eta\gamma_\sigma^-\rangle).$$
  This map is well defined
  since for different facet paths $\eta$ and $\eta'$ inducing the same coloring of the vertices of~$\eta$ we have
  $\langle\eta^-\eta'\rangle=1$.  Use the defining equation~(\ref{eq:gluing}) to conclude that the gluing in~$\oK'$ is
  equivalent to the gluing in~$\oK$.  The construction of~$\unfolding{K}'$ shows that $\unfolding{K}$ only depends on
  the combinatorial type of~$K$.
%
\end{proof}

We explicitly describe the equivalence relation arising on the disjoint union~$\oK'$ of admissibly colored facets.  Let
$(\sigma,b)$ and $(\tau,c)$ be colored facets and let $x$ be a point in the intersection~$\sigma\cap\tau$.  Let $\kappa$
be the unique simplex such that $x$ is contained its relative interior.  Then the point~$(x,b)\in(\sigma,b)$ is
identified with the point~$(x,c)\in(\tau,c)$ if and only if there exists a facet path $\gamma=(\sigma,\ldots,\tau)$ such
that all facets of~$\gamma$ lie in~$\st\kappa$ and $b = \langle\gamma\rangle c$.  In this case we say that the colorings
$b$ and~$c$ \emph{induce} each other in~$\st\kappa$.

Note that both definitions of the complete unfolding carry over to pseudo-simplicial complexes.

\begin{prp}
  The group of projectivities of the complete unfolding is trivial.
\end{prp}

\begin{proof}
  From the second construction of the complete unfolding each vertex of~$\oK'$ has a natural color.  The gluing process
  respects this coloring.  Therefore, the vertices of~$\unfolding{K}$ can be colored with $d+1$~colors.  By
  Proposition~\ref{prp:non-degenerate} there are no non-trivial projectivities in~$\unfolding{K}$.
\end{proof}

The following construction of the \emph{partial unfolding} is similar to the second definition of the complete
unfolding.  Let us consider the collection of all pairs $(\sigma,v)$, where $\sigma$ is a facet of $K$ and $v$ is a
vertex of $\sigma$.  As above we consider $(\sigma,v)$ as a geometric $d$-simplex affinely isomorphic to~$\sigma$.  Let
$\sigma$ and $\tau$ be neighbors. Then we glue $(\sigma,v)$ and $(\tau,w)$ along the common ridge of $\sigma$ and $\tau$
if $w=v\langle\sigma,\tau\rangle$.  As a result we obtain a pseudo-simplicial complex $\partunfold{K}$ which we call the
\emph{partial unfolding} of~$K$.  One can obtain the explicit description of the equivalence relation similar to that in
the case of the complete unfolding.  The partial unfolding map $r:\partunfold{K} \to K$ is induced by the affine
isomorphisms~$(\sigma,v)\to\sigma$.

In general, $\partunfold{K}$ is not connected.  We will denote by $\partunfold{K}_{(\sigma,v)}$ the component of
$\partunfold{K}$ containing the facet~$(\sigma,v)$.  It is immediate that $\partunfold{K}_{(\sigma,v)} =
\partunfold{K}_{(\tau,w)}$ if and only if there exists a facet path~$\gamma$ from $\sigma$ to~$\tau$ in~$K$ such that
$v\langle\gamma\rangle = w$.  In other words, the connected components of~$\partunfold{K}$ correspond to the orbits of
the action of~$\Pi(K,\sigma_0)$ on the set~$V(\sigma_0)$.

Observe the following properties of the unfoldings. The complete unfolding and each connected component of the partial
unfolding are strongly connected.  If $K$ is a pseudo-manifold (that is, each ridge of~$K$ is contained in exactly two
facets) then $\unfolding{K}$ and~$\partunfold{K}$ both are also pseudo-manifolds. If $K$ is orientable, then
$\unfolding{K}$ and~$\partunfold{K}$ both are also orientable.

The reader might have noted a similarity between the group of projectivities of a simplicial complex and the fundamental
group of a topological space. In the same spirit the complete unfolding is similar to the universal covering. This
analogy will become more evident in Section~\ref{sec:Odd}.

\subsection{Examples}\label{subsec:examples}

We give a few examples for the group of projectivities and the unfoldings.

\subsubsection{Graphs}\label{ex:graphs}

Graphs are precisely the $1$-dimensional simplicial complexes.  The group of projectivities of a graph~$\Gamma$ is
either trivial or it is isomorphic to~$S_2=\ZZ_2$ depending on whether $\Gamma$ is bipartite or not.  In the first case
the complete unfolding $\unfolding{\Gamma}$ is isomorphic to $\Gamma$, and the partial unfolding $\partunfold{\Gamma}$
consists of two copies of the graph $\Gamma$.  For a non-bipartite $\Gamma$ the complete unfolding coincides with the
partial unfolding, and it is a $2$-fold covering of~$\Gamma$.

\subsubsection{The starred triangle}\label{ex:starred}

\begin{figure}[htbp]
  \begin{center}
    \epsfig{file=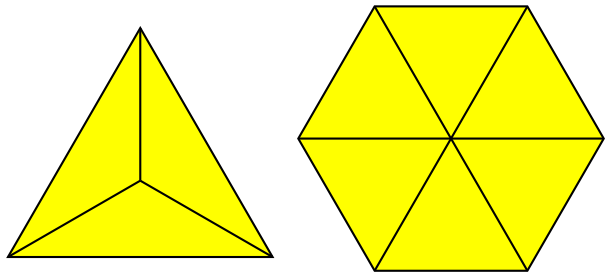,width=8cm}
    \caption{Starred triangle and its unfoldings.}
    \label{fig:simple-example}
  \end{center}
\end{figure}

Branching phenomena occur in dimension~$2$ and above.  Consider the cone~$\cT$ over the boundary of a triangle as in
Figure~\ref{fig:simple-example} (left).  We call this complex the \emph{starred triangle}.  Its group of projectivities
is generated by a transposition, which, for any base facet, exchanges the two vertices different from the apex.  The
complete unfolding is a $2$-fold covering with a unique branch point corresponding to the apex; the complex
$\unfolding{\cT}$ is isomorphic to the cone~$\cH$ over the boundary of a hexagon, see Figure~\ref{fig:simple-example}
(right).  The partial unfolding is the disjoint union of a copy of~$\cT$ (the unfolding with respect to the apex) and of
a copy of~$\cH$ (the unfolding with respect to any other vertex).

\subsubsection{The boundary of the $3$-simplex}

\begin{figure}[htbp]
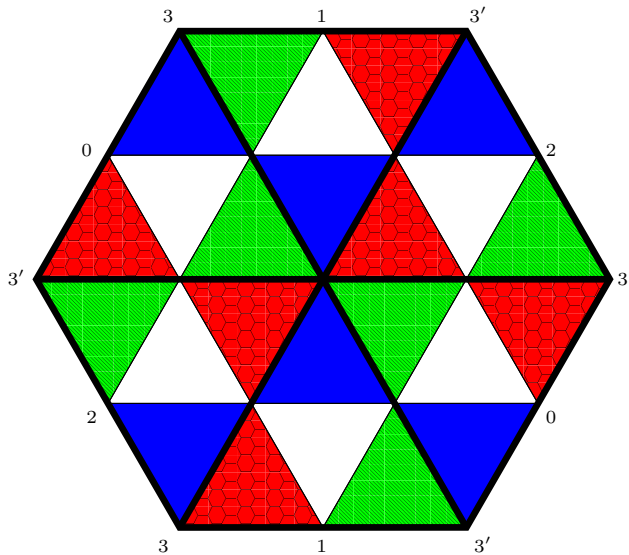

  \begin{center}
    \input unfolding-simplex.pstex_t
    \caption{Torus triangulation which arises as the complete unfolding of the boundary of the
      $3$-simplex.  The numbering of the vertices at the boundary of the hexagon indicate the identifications.  Facets
      of the same color belong to the same orbit under the action of~$\Pi(\partial\Delta^3)$.}
    \label{fig:unfolding-simplex}
  \end{center}
\end{figure}

The group of projectivities of the boundary complex $\partial\Delta^3$ of the $3$-dimensional simplex is the symmetric
group $S_3$.  The complete unfolding of $\partial\Delta^3$ is glued from $24$~triangles as follows.  We triangulate a
hexagon as shown in Figure~\ref{fig:unfolding-simplex} and then identify each pair of its opposite sides by translation.

Thereby, topologically $\unfolding{\partial\Delta^3}$ is a torus~$T^2$. The complete unfolding map
$\unfolding{\partial\Delta^3}\to\partial\Delta^3$ is a $6$-fold branched covering $T^2\to\Sph^2$ with $4$~branch points
on the sphere $\Sph^2$, the pre-image of each consists of $3$~points with branching index~$2$.

The partial unfolding $\partunfold{\partial\Delta^3}$ is the boundary of the tetrahedron with starred facets.  This
simplicial complex has $4$~vertices of degree~$3$ and $4$~vertices of degree~$6$.  Topologically the partial unfolding
is a $3$-fold branched covering of the $2$-sphere over itself with $4$~branch points.  The pre-image of each point
consists of one point with branching index~$2$ and of one point with branching index~$1$.

\subsubsection{A torus triangulation}\label{ex:torus}

\begin{figure}[htbp]
  \begin{center}
    \epsfig{file=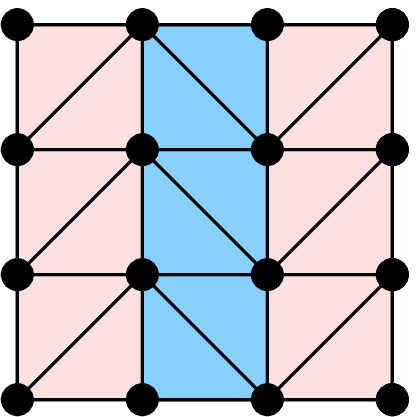,height=3cm}\qquad\epsfig{file=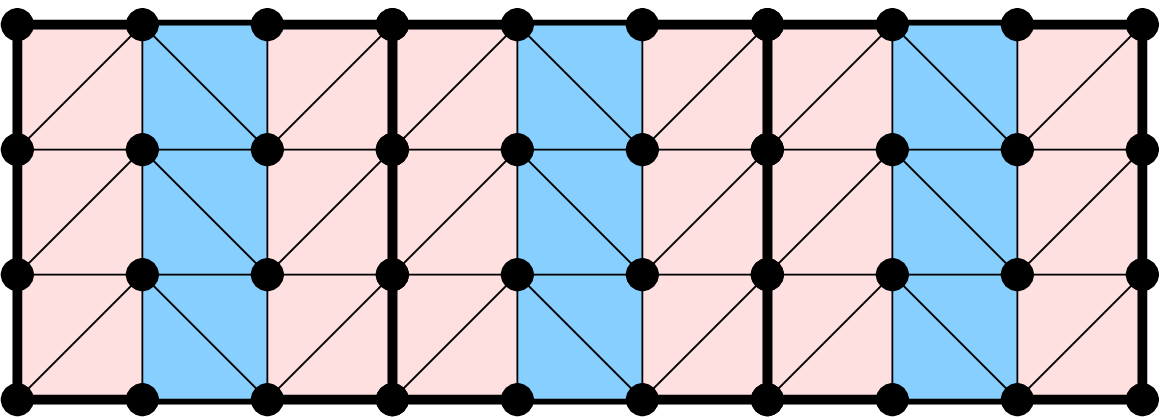,height=3cm}
    \caption{A torus triangulation and its unfolding.}
    \label{fig:anti-torus}
  \end{center}
\end{figure}

Branch points do not necessarily occur in high-dimensional unfoldings.  For example, consider a triangulation of the
$2$-torus as in Figure~\ref{fig:anti-torus} (left).  Its group of projectivities is cyclic of order~$3$.  The complete
and the partial unfoldings coincide.  Each of them is an unbranched $3$-fold covering, as shown in
Figure~\ref{fig:anti-torus} (right).

\subsubsection{A complex whose unfolding is not a simplicial complex}

\begin{figure}[htbp]
  \begin{center}
    \epsfig{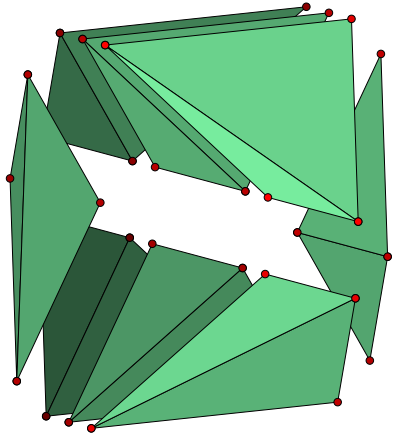}
    \caption{Explosion of a $3$-dimensional simplicial complex whose complete
      unfolding is not a simplicial complex.  The group of projectivities is trivial.  The complete
      unfolding is obtained from this complex by duplicating the middle horizontal edge. }
    \label{fig:unfolding-not-simplicial}
  \end{center}
\end{figure}

As mentioned above, an unfolding of a simplicial complex may not be a simplicial complex.  The first examples can be
found in dimension~$3$.  We outline the idea of a construction.  Consider two tetrahedra $\sigma$ and $\tau$ sharing a
common edge~$e$ with vertices $v$ and~$w$.  Let $b$ and $c$ be colorings of $\sigma$ and~$\tau$, respectively, such that
their restrictions to~$e$ coincide.  Then we can add further tetrahedra to $\sigma$ and~$\tau$ such that the following
holds: The colorings $b$ and~$c$ induce each other both within $\st v$ and within~$\st w$ but do not within $\st e$.
Then the colored facets $(\sigma,b)$ and~$(\tau,c)$ have the two vertices~$(v,b)=(v,c)$ and $(w,b)=(w,c)$ in common, but
no edge.  An example is shown on Figure~\ref{fig:unfolding-not-simplicial}. This simplicial complex is not locally
strongly connected in the sense of the definition given in Section~\ref{subsec:relation}.  However, the construction can
be modified to obtain a locally strongly connected example.

\subsection{Partial unfoldings of triangulations of the $3$-sphere include all $3$-manifolds}
\label{sec:formulation}

In this section we give a topological characterization of those branched coverings of $3$-manifolds which can be
obtained by unfoldings.  This implies our Main Result: Via the unfolding construction we obtain all closed orientable
$3$-manifolds from triangulations of the $3$-sphere.

\begin{thm}{\rm (Topological Characterization Theorem)}\label{thm:3dim}\\
  Let $N$ be a closed $3$-dimensional manifold, and let $f:M\to N$ be a branched covering with the following properties:
  \begin{enumerate}
  \item the number of sheets is less than or equal to~$4$;
  \item the branch set~$L\subset N$ is a knot or a link which is a boundary~$\mod 2$ in~$N$;
  \item the pre-image of each point in $L$ contains exactly one point of branching index~$2$; all other points in the
    pre-image are regular.
  \end{enumerate}
  Then there is a triangulation~$K$ of~$N$ such that $M$ is PL-homeomorphic to a component of the partial unfolding
  of~$K$ and $f$ is equivalent to the restriction of the partial unfolding map.
\end{thm}

The second property (the branch set~$L$ is a boundary~$\mod 2$ in $N$) means that the image of the fundamental cycle
of~$L$ under the natural homomorphism $H_1(L;\ZZ_2) \to H_1(N;\ZZ_2)$ is zero.

Each partial unfolding of a $3$-manifold has all the properties listed above.  In particular, by
Proposition~\ref{prp:PropertiesOfUnfolding} the branch set is always a boundary~$\mod 2$.

Now recall a theorem of Hilden and Montesinos (see~\cite{Hil} and~\cite{Mon}) which says that any closed orientable
$3$-manifold can be represented as a special kind of branched covering of the $3$-sphere.

\begin{thm}\label{thm:Hilden-Montesinos}
  Every closed orientable $3$-manifold~$M$ is a $3$-fold branched covering space of~$\Sph^3$ with a
  knot~$L$ as the branch set, such that the pre-image of each point of~$L$ consists of one point of
  branching index~$2$ and of one point of branching index~$1$.
\end{thm}

A glance at the conditions in Theorem \ref{thm:3dim} (together with the fact that $H_1(\Sph^3;\ZZ_2)=0$)
suffices to make the following conclusion.

\begin{thm}\label{thm:main-result}
  For each closed orientable $3$-manifold~$M$ there is a triangulation of the sphere~$\Sph^3$ such that one of the
  components of its partial unfolding is homeomorphic to~$M$.
\end{thm}

The Topological Characterization Theorem will be proven in Section~\ref{sec:3Unfolding}. Sections~\ref{sec:Odd}
to~\ref{sec:manifolds} contain the preliminaries which we need on the thorny path to the main results.

\section{Nice complexes and their unfoldings}\label{sec:Odd}

In this section we show how certain local properties of a pseudo-simplicial complex ensure a good behavior of its
unfoldings.  This should be seen in the context of coverings of topological spaces, where it is known that a satisfying
theory requires a variety of connectivity assumptions on the spaces involved.

\subsection{Relationship between the group of projectivities and the unfoldings}
\label{subsec:relation}

A simplicial complex~$K$ is called \emph{locally strongly connected} if the star of each face is strongly connected, see
also Mohar~\cite[p.~341]{Moh}.  In particular, this implies that $\partunfold{K}_{(\sigma,v)} =
\partunfold{K}_{(\tau,v)}$ for arbitrary facets $\sigma$ and~$\tau$ sharing a vertex~$v$.  Hence we denote this
component of the partial unfolding simply by $\partunfold K_v$.

A $d$-dimensional complex is \emph{balanced} if its vertices can be colored with $d+1$~colors so that there is no pair
of adjacent vertices with the same color.  In this case the coloring is unique up to renaming colors.  For combinatorial
properties of balanced complexes see Stanley~\cite[III.4]{Stanley}.

\begin{prp}\label{prp:Identity}
  Suppose that $K$ is a locally strongly connected simplicial complex.  Then the following are
  equivalent:
  \begin{enumerate}
  \item \label{item1} the group $\Pi(K)$ of projectivities is trivial;
  \item \label{item2} the complex $K$ is balanced;
  \item \label{item3} the complete unfolding map $p:\unfolding{K}\to K$ is a simplicial isomorphism;
  \item \label{item4} the restriction of the map $r:\partunfold{K}\to K$ to
    each component of $\partunfold{K}$ is a simplicial isomorphism.
  \end{enumerate}
\end{prp}

\begin{proof}
  The equivalence of the first two conditions was proved in~\cite[Proposition~6]{ProjSimplePolytopes}.
  
  Let us prove the equivalence of the third condition to the first one.  From the definition of the complete unfolding
  it is immediate that $|p^{-1}(x)|=|\Pi(K)|$ for any point~$x$ in the relative interior of any facet of~$K$.  Thus {\it
    (\ref{item3})\/} implies~{\it (\ref{item1})}.  On the other hand suppose that $\Pi(K)$ is trivial. Then
  $\unfolding{K}$ looks as follows. Take the disjoint union~$\Sigma(K)$ of all facets of $K$ and, for each pair of
  facets which are neighbors in $K$, glue them along the common ridge. Note that the complex~$K$ can also be obtained
  from $\Sigma(K)$ in a similar way, with the only difference that gluings must be performed not only along the ridges,
  but along faces of all dimensions.  We must show that the two equivalence relations on~$\Sigma(K)$ defined in this way
  are the same.  Obviously, the second equivalence relation is stronger or equal than the first one.  Conversely, let
  $\kappa$ be a face of $K$, and let $\kappa\subset\sigma\cap\tau$, where $\sigma$ and~$\tau$ are facets.  Since $K$ is
  locally strongly connected, there is a facet path from $\sigma$ to~$\tau$ such that the face $\kappa$ lies in all
  facets of this path.  This implies that in the complete unfolding $\unfolding{K}$ the facets $\sigma$ and $\tau$ are
  glued along the face $\kappa$.
  
  Now proceed to the fourth condition.  As it was already mentioned, the components of the partial unfolding are in the
  one-to-one correspondence with the orbits of the $\Pi(K,\sigma_0)$-action on the set~$V(\sigma_0)$.  Besides,
  $|r^{-1}(x)|=d+1$ for a point~$x$ in the relative interior of any facet of $K$.  Thus, if each component of
  $\partunfold{K}$ is mapped to~$K$ isomorphically, there must be exactly $d+1$~orbits and $\Pi(K)$ is trivial.  This
  shows that {\it (\ref{item1})\/} follows from~{\it (\ref{item4})\/}.  Finally, we prove that {\it (\ref{item2})\/}
  implies~{\it (\ref{item4})}.  Suppose that the vertices of the complex~$K$ are $(d+1)$-colored.  Then any component
  of~$\partunfold{K}$ is composed from the facets~$(\sigma,v)$, where $v$ ranges over the set of vertices of a fixed
  color.  Note that for any facet~$\sigma$ of~$K$ there is a unique vertex~$v$ of a given color.  The rest of the proof
  is similar to the argument given for the implication {\it (\ref{item1})} $\Rightarrow$ {\it (\ref{item3})}.
\end{proof}

\begin{figure}[htbp] 
  \begin{center}
    \epsfig{file=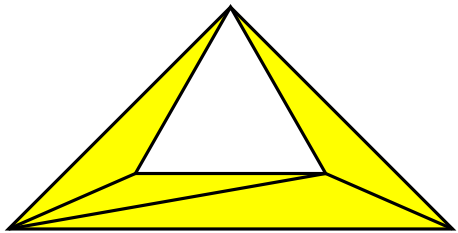,width=6cm}
    \caption{Strongly connected but not locally strongly connected complex.}
    \label{fig:NotLoc}
  \end{center}
\end{figure}

The $2$-dimensional complex in Figure~\ref{fig:NotLoc} is strongly connected but not locally strongly connected: The
star of the top vertex is not strongly connected.  One can see that the group of projectivities is trivial, although the
complex is neither balanced nor isomorphic to its complete unfolding (the vertex at the top is duplicated in the
unfolding).

A characterization of local strong connectivity is given by the following.

\begin{lem}
  A simplicial complex~$K$ is locally strongly connected if and only if for each face~$\kappa$ of~$K$ with
  $\codim\kappa>1$ the link $\lk\kappa$ is connected.
\end{lem}

\begin{proof}
  For any face~$\kappa$ with $\codim\kappa>1$ the star $\st\kappa$ is strongly connected if and only if the link
  $\lk\kappa$ is strongly connected.  Besides, a strongly connected complex is connected.  This proves that the above
  condition on links holds for any locally strongly connected complexes.
  
  To prove the converse implication suppose that $\kappa$ is a face in $K$ which is maximal (by inclusion) among the
  faces whose star is not strongly connected.  Clearly, $\codim\kappa>1$.  Put $L=\lk\kappa$.  The complex $L$ is pure,
  connected, but not strongly connected.  It easily follows that $L$ is not locally strongly connected.  Let $\lambda\in
  L$ be such that $\st_L\lambda$ is not strongly connected.  Note that
  \begin{equation} \label{eq:LinkLink}
    \lk_L\lambda=\lk_K(\kappa*\lambda).
  \end{equation}
  Since $\codim_K (\kappa*\lambda)=\codim_L\lambda>1$, we have that $\st(\kappa*\lambda)$ is not strongly connected. But
  this contradicts the assumption that $\kappa$ is a maximal face with this property.
\end{proof}

For homotopy properties of locally strongly connected complexes see Section~\ref{subsec:homotopy}.

\subsection{Relationship between $\Pi(K)$ and $\pi_1(K)$}

The link $\lk\kappa$ of a codimension-$2$-face $\kappa\in K$ is a graph which is connected provided that $K$ is locally
strongly connected.  Whenever this graph is bipartite, $\kappa$ is called an \emph{even} face, otherwise $\kappa$ is
called \emph{odd}.  The collection of all odd codimension-$2$-faces together with all their proper faces is called the
\emph{odd subcomplex} of $K$ and denoted by $\odd{K}$.  The odd subcomplex is pure, and it has codimension~$2$ or it is
empty.

For each face in~$K$ there is a natural correspondence between the facets in the star and the facets in the link.  This
correspondence extends to facet paths and thus to projectivities.  Hence we obtain a canonical isomorphism between the
groups of projectivities.  In particular, for a codimension-$2$-face $\kappa$, the group
$\Pi(\st\kappa)\cong\Pi(\lk\kappa)$ vanishes if and only if $\kappa$ is even; see the Example~\ref{ex:graphs}.  Thus, in
order to have $\Pi(K)=0$ it is clearly necessary that $\odd{K}=\emptyset$.

As a side remark, without a proof, we state the following.

\begin{prp}
  Suppose that $K_{\mathrm{odd}}$ is a locally strongly connected pseudo-manifold (equivalently, each $1$-dimensional
  link in $\odd{K}$ is a circle). Then $\odd{(\odd{K})}=\emptyset$.
\end{prp}

Observe that each homology manifold is a locally strongly connected pseudo-manifold.  In particular, each PL-manifold is
of this type.  For an introduction to homology manifolds see Munkres~\cite[\S63]{Mun}.

We call the simplicial complex $K$ \emph{locally strongly simply connected} if for each face $\kappa$ with
$\codim\kappa>2$ the link of $\kappa$ is simply connected.  Further, we call a complex \emph{nice} if it is locally
strongly connected and locally strongly simply connected.  Nice complexes can be seen as combinatorial analogues to
`sufficiently connected' topological spaces in the theory of coverings.  Observe that the class of nice complexes
contains all combinatorial manifolds as well as all graphs.

The key result~\cite[Theorem~8]{ProjSimplePolytopes} indicates how the (reduced) group of projectivities of a
combinatorial manifold is generated by special projectivities.  Here we are after a generalization to arbitrary nice
complexes.  For the proof of~\cite[Theorem~8]{ProjSimplePolytopes} it was convenient to work in the dual cell complex of
a combinatorial manifold.  While it is possible to define the \emph{dual block complex} for an arbitrary simplicial
complex (see Appendix~\ref{subsec:homotopy}), the resulting blocks do not have a good topological structure in general.
Our niceness condition is devised to make sure that the relative homotopy type of each dual block with respect to its
boundary is well behaved.

From now on all our complexes are supposed to be nice.  Fix a facet $\sigma_0$ in~$K$.

Let $\kappa$ be a codimension-$2$-face, $\sigma$ a facet containing~$\kappa$, and $g$ a path from $\sigma_0$
to~$\sigma$.  Since $\st\kappa$ is always simply connected we infer that the path $glg^-$ is null-homotopic for any
facet loop~$l$ in~$\st\kappa$ based at~$\sigma$.  Thus we have $\langle glg^-\rangle\in\Pi_0(K,\sigma_0)$.  We call such
a projectivity a projectivity \emph{around}~$\kappa$.  A projectivity around a codimension-$2$-face~$\kappa$ is either a
transposition or the identity map, depending on $\kappa$ being odd or even.

\begin{thm}
  The reduced group of projectivities $\Pi_0(K,\sigma_0)$ of a nice complex~$K$ is generated by the projectivities
  around the odd codimension-$2$-faces.  In particular, $\Pi_0(K,\sigma_0)$ is generated by transpositions.
\end{thm}

\begin{proof}
  Let $\gamma$ be a null-homotopic facet path in~$K$.  This yields a closed PL path $\overline\gamma:\Sph^1\to K$.
  
  Now, let $K^*_{(m)}$ denote the induced subcomplex of the barycentric subdivision $b(K)$ which is generated by the
  vertex set $\SetOf{\hat\tau}{\codim\tau\le m}$.  This is the \emph{dual $m$-skeleton} of~$K$, see also Appendix.
  Recalling the definition of the map~$\overline\gamma$ it is obvious that $\overline\gamma(\Sph^1)\subset K^*_{(1)}$,
  that is, $\overline\gamma$ is a closed path in the dual graph of~$K$.
  
  Applying Proposition~\ref{prp:BlockApprox} to the null-homotopic map~$\overline\gamma$ we find a map $g:\DD^2\to
  K^*_{(2)}$ such that the restriction of~$g$ to~$\Sph^1$ equals~$\overline\gamma$.  Here $\DD^2$ denotes the closed
  unit disk in~$\RR^2$ with $\partial\DD^2=\Sph^1$.  Finally, due to the Relative Simplicial Approximation Theorem there
  exists a PL map $h:\DD^2\to K^*_{(2)}$ which coincides with $g$ on $\Sph^1$.
  
  Consider the subpolyhedron
  $$C=\bigcup\SetOf{h^{-1}(\hat\kappa)}{\text{$\kappa$ face of codimension~$2$}}$$
  of the disk~$\DD^2$.  Let
  $C_1,\ldots,C_s$ be the connected components of $C$.  Subdivide the disk $\DD^2$ into polyhedral disks
  $D_1,\ldots,D_s$ using PL paths starting from the base point of $\Sph^1$ such that $C_i\subset\int D_i$ for
  $i=1,\ldots,s$.  This is shown in Figure~\ref{fig:subdivide}.  For each $i$ let $\phi_i$ be a simple closed path
  running along the boundary of the disk $D_i$ co-oriented with the boundary circle of $\DD^2$.  Clearly, the loop
  $\overline\gamma$ is homotopic to the product $\prod_{i=1}^s h\circ\phi_i$.

  \begin{figure}[htbp]
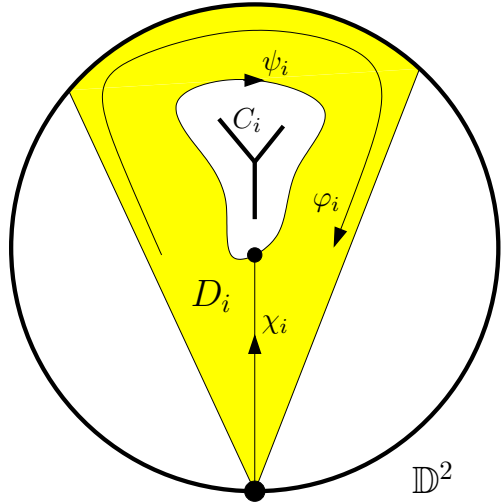
 
    \begin{center}
      \input subdivide.pstex_t
      \caption{Subdvision of the disk~$\DD^2$ and homotopy of paths.}
      \label{fig:subdivide}
    \end{center}
  \end{figure}
  
  We need an interpretation of continuous paths as facet paths. Suppose that a path $\xi:[0,1]\to K$ has the following
  properties:
  \begin{enumerate}
  \item $\xi([0,1])\cap K^{(d-2)}=\emptyset$;
  \item $\xi(0),\, \xi(1)\notin K^{(d-1)}$;
  \item the set $\xi^{-1}(K^{(d-1)})\subset(0,1)$ has a finite number of connected components.
  \end{enumerate}
  Note that the last condition holds for any PL path.  Under these assumptions we can form a facet path from the
  sequence of the facets encountered by $\xi$ on its way.  This will be called the \emph{discretization} of the path
  $\xi$.  For example, $\gamma$ is the discretization of the PL path~$\overline\gamma$.  To simplify notation we usually
  identify a suitable PL path with its discretization, in particular, we write $\gamma=\overline\gamma$.
  
  Suppose we have a homotopy $\xi_t:[0,1]\to K$ with fixed endpoints, where $t\in[0,1]$, such that each path $\xi_t$
  satisfies the above conditions.  Then it follows that the projectivities $\langle\xi_0\rangle$ and
  $\langle\xi_1\rangle$ coincide.
  
  Each of the paths $\gamma$, $h\circ\phi_i$ satisfies the conditions mentioned above.  Thus the statement just
  formulated implies the equality $$\langle\gamma\rangle = \prod_{i=1}^s \langle h\circ\phi_i \rangle.$$
  
  Further, each $\phi_i$ is PL homotopic inside $D_i$ to a path of the form $\chi_i\psi_i\chi_i^{-1}$, where the path
  $\psi_i$ is sufficiently close to $C_i$. Namely we assume that the path $h\circ\psi_i$ lies in the barycentric star of
  $\hat{\kappa_i}$, where $\hat{\kappa_i}=h(C_i)$. Then all the facets in the discretization of $h\circ\psi_i$ are in
  $\st\kappa_i$ and since
  $$
  \langle h\circ\phi_i\rangle = \langle h\circ\chi_i\rangle\langle h\circ\psi_i\rangle\langle
  h\circ\chi_i^{-1}\rangle
  $$
  we see that $\langle\gamma\rangle$ equals to a product of projectivities around codimension-2-faces.
\end{proof}

\begin{cor}\label{cor:vanish}
  If $\pi_1(K)$ is trivial and $\odd{K}=\emptyset$ then $\Pi(K)$ is trivial.
\end{cor}

The proof is straightforward.

\begin{cor}\label{cor:Odd}
  The odd subcomplex $\odd{K}$ coincides with the collection of faces~$\kappa$ such that $\Pi(\st\kappa)$ is
  non-trivial.
\end{cor}

\begin{proof}
  It is easy to see that this collection is indeed a simplicial complex.  Suppose that $\Pi(\st\kappa)$ is non-trivial.
  Then it is sufficient to prove that $\kappa$ is contained in a codimension-$2$-face with a non-bipartite link.
  Clearly, $\codim\kappa\ge 2$.  In the case of equality $\lk\kappa$ is not a bipartite graph.  Let $\codim\kappa>2$.
  Since $K$ is nice $\pi_1(\lk\kappa)$ is trivial and by Corollary~\ref{cor:vanish} the subcomplex $\odd{(\lk\kappa)}$
  is non-empty.  Thus there exists a face $\lambda\in\lk\kappa$ such that $\lk_{\lk\kappa}\lambda$ is a non-bipartite
  graph.  Since $\lk_{\lk\kappa}\lambda=\lk_K(\kappa*\lambda)$, the join $\kappa*\lambda$ is the desired odd
  codimension-$2$-face of~$K$.
\end{proof}

\subsection{Local behavior of the unfoldings}

Since the identifications which define the unfoldings have a local character (they are completely defined by the
structure of~$\st\kappa$), the unfoldings of some neighborhood of~$\kappa$ can be described in terms of~$\st\kappa$.
Our aim now is to give such a description.  Let us introduce the following notation.  The union of the relative
interiors of the faces containing a given face~$\kappa$ is denoted by~$\ost\kappa$. It is called the \emph{star of the
  interior} of~$\kappa$:
$$
\ost\kappa=\bigcup_{\sigma\supseteq\kappa}\relint\sigma.
$$

In other words $\ost\kappa$ is the complement to the union of those faces of $K$ which do not contain $\kappa$. Thus it
is an open subset of $K$. We have the identity
$$
\ost\kappa = (\kappa*\lk\kappa)\setminus(\partial\kappa*\lk\kappa).
$$

\begin{lem}\label{lem:StarOfInterior}
  Let $K$ be a nice complex with complete and partial unfoldings $p:\unfolding{K}\to K$ and $r:\partunfold{K}\to K$,
  respectively.  For each face~$\kappa$ of~$K$ the following holds:
  \begin{enumerate}
  \item Each component of the pre-image $p^{-1}(\ost\kappa)$ is homeomorphic to the space
    $$
    (\kappa*\unfolding{\lk\kappa})\setminus(\partial\kappa*\unfolding{\lk\kappa}).
    $$
    Moreover, the homeomorphism between a component of $p^{-1}(\ost\kappa)$ and the above space can be chosen such
    that the projection onto~$\ost\kappa$ is induced by the natural map $\unfolding{\lk\kappa}\to\lk\kappa$.
  \item Each component of the pre-image $r^{-1}(\ost\kappa)$ is either homeomorphic to $\ost\kappa$
    or to the space
    $$
    (\kappa*\partunfold{\lk\kappa}) \setminus(\partial\kappa*\partunfold{\lk\kappa}),
    $$
    where the projection of the component onto $\ost\kappa$ is induced either by identity map or by the natural map
    $\partunfold{\lk\kappa}\to\lk\kappa$.
  \end{enumerate}
\end{lem}

\begin{proof}
  To specify a component~$C$ of $p^{-1}(\ost\kappa)$ it is sufficient to pick a facet~$(\sigma,b)$ of $\unfolding{K}$
  having a non-empty intersection with that component.  Here $b$ is an admissible coloring of the facet~$\sigma$ of~$K$.
  Clearly, then $\sigma\in\st\kappa$.  For another facet $(\tau,c)$ the intersection
  $(\sigma\cap\ost\kappa,b)\cap(\tau\cap\ost\kappa,c)$ is either empty or equals $(\relint\kappa,b)=(\relint\kappa,c)$.
  The latter is the case if and only if the colorings $b$ and~$c$ induce each other in~$\st\kappa$.  Hence the
  component~$C$ is fromed by the blocks $(\tau\cap\ost\kappa,c)$ such that $c$ is induced by~$b$ in~$\st\kappa$.
  
  There is a natural bijection between the facet paths in~$\st\kappa$ and those in~$\lk\kappa.$ Since this bijection
  respects the inducing of colorings, the identifications between the blocks functorially arise from the identifications
  at the construction of the space $\unfolding{\lk\kappa}$.  The resulting space is homeomorphic to
  $(\kappa*\unfolding{\lk\kappa})\setminus(\partial\kappa*\unfolding{\lk\kappa})$.
  
  The proof of the second part of the lemma is similar.  We specify a component of $r^{-1}(\ost\kappa)$ by picking a
  facet $(\sigma,v)$, where $\sigma\in\st\kappa$.  Now we have two cases.  If $v\in\kappa$, then the component is
  homeomorphic to~$\ost\kappa$.  Otherwise, if $v\notin\kappa$, then arguing as in the previous paragraph we get the
  space functorially arising from $\unfolding{\lk\kappa}$.
\end{proof}

Now we are ready to prove the following key result.

\begin{thm}\label{thm:UnfoldingIsBranched}
  The restriction of the complete unfolding of~$K$ to the pre-image of the complement of the odd subcomplex is a
  covering.  The same is true for each component of the partial unfolding.
\end{thm}

\begin{proof}
  We must prove that the map $(\kappa*\unfolding{\lk\kappa}) \setminus(\partial\kappa*\unfolding{\lk\kappa})\to
  \ost\kappa$ is a homeomorphism for any $\kappa\notin\odd{K}$. It will follow from the statement that the unfolding
  $\unfolding{L}\to L$, where $L=\lk\kappa$, is a simplicial isomorphism. This is really the case due to
  Proposition~\ref{prp:Identity} and Corollary~\ref{cor:Odd}.  However, in order to apply these propositions we must
  make sure that $L$ is nice. But indeed, $L$ is pure since $K$ is pure, $L$ is strongly connected since $K$ is locally
  strongly connected and, finally, due to (\ref{eq:LinkLink}) the local strong connectedness of $L$ follows from that of
  the complex $K$.
  
  It should be also shown that the space $K\setminus\odd{K}$ and its pre-images in the complete unfolding and in each
  component of the partial unfolding are connected.  This easily follows from the strong connectedness of~$K$ and of the
  unfoldings.
\end{proof}

\section{The unfoldings are branched coverings}\label{sec:BrCov}

In this section we show that the unfoldings of a nice simplicial complex~$K$ are branched coverings in the sense of
Fox~\cite{Fox}.  Moreover, it turns out that the complete unfolding can be related to the partial unfolding in a purely
topological way.  That is to say, the relationship does not depend on the combinatorial structure of~$K$.  Below we
assume that $\unfolding{K}$ and $\partunfold{K}$ are simplicial complexes (not only pseudo-complexes).  We may do so
since we can subdivide~$K$ as described in the Appendix~\ref{subsec:anti-prismatic}.  The key property of this
subdivision is that the unfoldings of the subdivided complex turn out to be PL equivalent to the unfoldings of~$K$.

\subsection{Branched coverings}\label{subsec:BranchedCoverings}

The topological concept of a branched covering was formulated by Fox~\cite{Fox}.  In his approach branched coverings are
derived from (unbranched) coverings.  As usual in this section the notion of a covering is restricted to the case where
the covering space is connected.

Let $h:X\to Z$ be a continuous map with the following properties.  Firstly, the restriction $h:X\to h(X)$ is a covering.
Secondly, the image~$h(X)$ is a dense subset of~$Z$.  And, thirdly, $h(X)$ is \emph{locally connected in}~$Z$, that is,
in each neighborhood~$U$ of each point $z\in Z$ there is an open set $V\ni z$ such that the intersection $V\cap h(X)$ is
connected.  In this situtation Fox constructs a \emph{completion} of $h$ which is a surjective map $g:Y\to Z$ with
$Y\supseteq X$ and $g|_X=h$.  Any two completions $g_i:Y_i\to Z$ for $i=1,2$ are equivalent in the sense that there
exists a homeomorphism $\phi:Y_1\to Y_2$ satisfying $g_2\circ\phi=g_1$ and $\phi|_X=\id$.  By definition, any map
$g:Y\to Z$ obtained in this way is called a \emph{branched covering}.  It may happen that $g$ again is a covering.  In
this way a covering is a special case of a branched covering.

If $g:Y\to Z$ is an arbitrary map, then let $Z_{\rm ord}$ denote the unique maximal subset of~$Z$ such that the
corresponding restriction of~$g$ is a covering over $Z_{\rm ord}$.  Put $Z_{\rm sing}=Z\setminus Z_{\rm ord}$.  If $g$
is a branched covering then $Z_{\rm sing}$ is called the \emph{singular} or \emph{branch set} of~$g$.  And the
restriction $g_{\rm ord}:g^{-1}(Z_{\rm ord})\to Z_{\rm ord}$ is called the covering \emph{associated} with~$g$.  It is
the maximal covering whose completion is equivalent to~$g$.  For a simplicial map $f:J\to K$, Fox presents necessary and
sufficient conditions for $f$ to be a branched covering.  Note that the singular set $K_{\rm sing}$ of a simplicial map
is a subcomplex of~$K$.

\begin{prp}[Fox~\cite{Fox}, p.~251] \label{prp:SimplicBranched}
  A simplicial map $f:J\to K$ is a branched covering if and only if the following conditions hold:
\begin{enumerate}
\item $f$ is non-degenerate, that is, $f$ maps each simplex onto a simplex of the same dimension;
\item for each face $\tau$ of~$K_{\rm sing}$ the space $K_{\rm
    ord}\cap\st_K\sigma$ is non-empty and connected;
\item $f^{-1}(K_{\rm ord})$ is connected;
\item for each face~$\sigma$ in~$f^{-1}(K_{\rm sing})$ the space $f^{-1}(K_{\rm
    ord})\cap\st_J\sigma$ is non-empty and connected.
\end{enumerate}
\end{prp}

For simplicial maps between triangulated manifolds these conditions are equivalent to the classical one: $K_{\rm sing}$
has codimension at least~$2$.  Moreover, the same reformulation holds for a wider class of simplicial complexes, see
also Mohar~\cite[p.~341]{Moh}.

\begin{prp}
  If the simplicial complexes~$J$ and $K$ are pure, strongly connected, and locally strongly connected, then the
  conditions in the previous Proposition are equivalent to the inequality
  \begin{equation}\label{eq:CodimSing}
    \codim K_{\rm sing}\geq 2.
  \end{equation}
\end{prp}

The proof is left to the reader.

Now, by Theorem~\ref{thm:UnfoldingIsBranched}, both the complete unfolding of~$K$ and each component of the partial
unfolding are branched coverings.  The branch set of the complete unfolding is the odd subcomplex~$\odd{K}$.

The following is a specialization of a definition given by Fox~\cite{Fox}.  Let $f:J\to K$ be a branched covering, where
$J$ and $K$ both are pure, strongly connected, and locally strongly connected simplicial complexes.  Let $\sigma$ be any
face of $f^{-1}(K_{\rm sing})$.  Denote $\tau=f(\sigma)$ and let $O$ be a connected component of $f^{-1}(\ost\tau)$
which contains $\relint\sigma$.  Then the number of sheets of the branched covering $f|_O:O\to\ost\sigma$ is called the
\emph{index of branching} at the face $\sigma$.  We write $\ind_f\sigma$.  Additionally, the index of branching at an
arbitrary point $x\in f^{-1}(K_{\rm sing})$ is defined to be $\ind_f\sigma$ where $\sigma$ is the unique face with
$x\in\relint\sigma$.

\subsection{The complete unfolding is the regularization of the partial unfolding}

In this section we allow disconnected covering spaces. Thus we can view the partial unfolding as a branched covering as
well.

Since a branched covering has a uniquely determined associated covering, we can carry over some familiar concepts of the
theory of coverings to the more general branched case.  Thus, for example, $|\Pi(K)|$ is the number of sheets of the
branched covering $p:\unfolding{K}\to K$.  Similarly, for $r:\partunfold{K}\to K$ the number of sheets is equal
to~$d+1$.

Moreover, the notion of a covering transformation naturally generalizes to branched coverings: A \emph{covering
  transformation}, or \emph{Deck transformation}, is a homeomorphism of the (branched) covering space which commutes
with the projection.  Each covering transformation of a branched covering arises as the unique extension of a covering
transformation of the associated covering.  Hence a branched covering $f$ has the same group~$\frD(f)$ of covering
transformations as its associated covering.

Let us recall some facts about coverings.  Throughout the following let $f:X\to Y$ be a covering, where $X$ and~$Y$
satisfy certain connectivity properties, see Bredon~\cite[III.3.1 and III.8.3]{Bred}.  We choose a base point $x_0\in
X$, and we put $y_0=f(x_0)$.  Then $f$ induces a monomorphism $f_*:\pi_1(X,x_0)\to\pi_1(Y,y_0)$.  The image is called
the \emph{characteristic subgroup} for~$f$.  The equivalence classes of coverings over~$Y$ are in one-to-one
correspondence with the conjugacy classes of subgroups of~$\pi_1(Y,y_0)$.  The covering~$f$ is \emph{regular} if its
characteristic subgroup is a normal subgroup.  For any subgroup $H$ of~$G$ let $\core_G(H)=\bigcap_{g\in G}g^{-1}Hg$ be
the \emph{core} of $H$ in~$G$.  This is the largest normal subgroup of~$G$ which is contained in~$H$.  The
\emph{regularization} of the covering~$f$ is the (regular) covering corresponding to the core of the characteristic
subgroup of~$f$ in~$\pi_1(Y,y_0)$.  We call a branched covering \emph{regular} if its associated covering is regular.
Similarly, one branched covering is a \emph{regularization} of another if the same holds for their associated coverings.

In the next two paragraphs we describe two tools for classification of coverings over fixed space $Y$. The first one
enumerates all coverings, the second one --- only regular coverings.

Let $\gamma:[0,1]\to Y$ be a path in~$Y$, and let $x\in X$ with $f(x)=\gamma(0)$.  Then there is a unique
path~$L(\gamma,x)$ in~$X$ with $f\circ L(\gamma,x)=\gamma$ starting at~$x$.  We denote its endpoint by~$e(\gamma,x)$.
For any closed path~$\gamma$ starting at~$y_0$ and $x_i\in f^{-1}(y_0)$ the point~$e(\gamma,x_i)$ is contained
in~$f^{-1}(y_0)$.  As this only depends on the homotopy class of~$\gamma$, we obtain an action of the fundamental
group~$\pi_1(Y,y_0)$ on the set~$f^{-1}(y_0)$.  This defines the \emph{monodromy homomorphism}
$$\frm_f:\pi_1(Y,y_0)\to\Sym(f^{-1}(y_0)).$$
The image $\frM(f)=\frm_f(\pi_1(Y,y_0))$ is called the \emph{monodromy
  group} of~$f$; see Seifert and Threlfall~\cite[\S58]{ST}.  Conversely, for any homomorphism $\frm:\pi_1(Y)\to S_n$,
where $S_n$ denotes the symmetric group of degree~$n$, there is an (up to equivalence unique but not necessarily
connected) $n$-fold covering~$f$ such that $\frm_f=\frm$.  Moreover, conjugation in~$S_n$ does not change the
equivalence class of~$f$.

Suppose that $f$ is regular.  Then there is an epimorphism $\frd_f:\pi_1(Y,y_0)\to\frD(f)$ defined as follows.  For any
point $x\in X$ let $\eta'$ be a path in~$X$ from $x_0$ to $x$, and let $\eta=f\circ\eta'$ be its projection to~$Y$.  We
have $\eta(0)=y_0$ and $\eta(1)=f(x)$.  Now the covering transformation $\frd_f[\gamma]$ maps the point~$x$
to~$e(\eta^{-1}\gamma\eta,x)$.  This definition does not depend on the choice of~$\eta'$.  Conversely, if a group $G$
and an epimorphism $\frd:\pi_1(Y)\to G$ are given, then there exists an (up to equivalence unique) regular covering~$f$
over~$Y$ with $\frd_f=\frd$; the group~$G$ coincides with the group of covering transformations.  We call $\frd_f$ the
\emph{characteristic homomorphism} of $f$.  Observe that $\ker\frd_f$ is the characteristic subgroup for~$f$.

One can prove that the regularization of the covering $f$ is the regular covering whose characteristic homomorphism is
the monodromy homomorphism of~$f$.


\begin{thm}
  The complete unfolding is the regularization of the partial unfolding.
\end{thm}

\begin{proof}
  The group $\Pi(K)$ acts on the set $\Sigma(K)\times\Pi(K)$.  Due to the equation~(\ref{eq:gluing}) this action
  descends to an action on the complete unfolding.  Clearly, each fiber is invariant under this action.  For any
  point~$a$ in the relative interior of some facet of~$K$ consider the fiber $p^{-1}(a)$, where $p$ is the complete
  unfolding map.  The action of~$\Pi(K)$ on this fiber is equivalent to the action of~$\Pi(K)$ on itself by
  multiplication on the right.  In particular, this action is transitive.  On the other hand, the covering
  transformation group acts freely on each non-singular fiber.  We conclude that $\frD(p)=\Pi(K)$.
  
  Let $\sigma$ be the facet which contains~$a$.  For the partial unfolding map~$r$, the fiber~$r^{-1}(a)$ is the set
  $\{a\}\times V(\sigma)$.  It follows that the action of~$\Pi(K,\sigma)$ on the set~$V(\sigma)$ is the monodromy
  action.
  
  We denote the homomorphism $\frd_p=\frm_r:\pi_1(K\setminus\odd{K},a)\to\Pi(K,\sigma)$ by $\frh_K$.  It has the
  following form:
  \begin{equation}
  \frh_K[\gamma]=\langle\gamma\rangle,\label{eq:frh}
  \end{equation}
  where $\gamma$ on the right is an arbitrary facet path, and $\gamma$ on the left is the corresponding path in the dual
  graph.  Since, for a nice complex, the paths in the dual graph generate the fundamental group (see Proposition
  \ref{prp:PathApprox} in the Appendix), the homomorphism $\frh_K$ is determined by this equality.
\end{proof}

The above theorem shows that topologically the complete unfolding can be derived from the topological type of the
partial unfolding.

The following result will not be used in the sequel, and we mention it without proof which is straightforward. It says
that the complete unfolding is a composition of partial unfoldings and thus provides another connection between these
objects.

First make some conventions.  Choose a facet $\sigma_0$ in a simplicial complex~$K$ (which need not to be nice) and an
arbitrary ordering $(v_0,\ldots,v_d)$ of the vertices of $\sigma$. In the partial unfolding of~$K$ take the component
$\partunfold K_{(\sigma,v_0)}$ and denote it by $K_1$.  Denote the facet $(\sigma,v_0)$ of the complex $K_1$ by $\sigma_1$.
The vertices of $\sigma_1$ are in a natural correspondence with those of $\sigma_0$ and, by abuse of notation, we denote
them in the same way.  Now in the partial unfolding of~$K_1$ we take the component containing the facet
$(\sigma_1,v_1)$. Proceeding in this manner we obtain a sequence of pseudo-complexes
\begin{equation}
\label{eqn:tower}
K=K_0\stackrel{r_0}{\leftarrow} K_1\stackrel{r_1}{\leftarrow}\cdots
\stackrel{r_d}{\leftarrow} K_{d+1},
\end{equation}
where $K_{i+1}=\partunfold{(K_i)}_{(\sigma_i,v_i)}, \sigma_{i+1} = (\sigma_i,v_i)$, for $i=0,\ldots,d$.  Here $r_i$ stands
for a restriction of the corresponding partial unfolding map.

\begin{thm}\label{thm:composition}
  The composition map $K_{d+1}\to K$ is simplicially equivalent to the complete unfolding, that is, there exists a
  simplicial isomorphism between pseudo-complexes $K_{d+1}$ and $\unfolding{K}$ such that the corresponding diagram
  commutes.
\end{thm}

\section{Unfoldings of PL-manifolds}
\label{sec:manifolds}

An important class of nice complexes are triangulations of PL-manifolds.  Throughout, we tacitly assume such
triangulations to be compatible with the fixed PL-structure.  A fixed triangulation of a PL-manifold is occasionally
called a combinatorial manifold, e.g., see Glaser~\cite{Gla}.  A \emph{subpolyhedron} of a PL-manifold is the subspace
corresponding to a subcomplex of some triangulation.

Suppose we have a branched covering $f:M\to N$, where $N$ is a PL manifold and $N_{\rm sing}$ is a subpolyhedron of~$N$.
Then the map $f$ induces a polyhedral structure on~$M$.  The following problem arises: Give necessary and sufficient
conditions on the associated covering $f_{\rm ord}$ such that $M$ is a PL-manifold.

Fox~\cite{Fox} shows that the polyhedron~$M$ is a PL-manifold provided that $N_{\rm sing}$ is a locally flat
codimension-$2$-submanifold and the index of branching is everywhere finite.  A submanifold~$L\subset N$ is
\emph{locally flat} if for each point $x\in L$ there is a neighborhood~$U\subset N$ such that the pair $(U,U\cap L)$ is
PL-homeomorphic to the pair~$(\DD^m,\DD^n)$.  A $1$-dimensional submanifold is always locally flat.  For some necessary
conditions see Hemmingsen~\cite{Hem}.

Here we are concerned with the question which branched coverings of PL-manifolds with a locally flat branch set arise as
unfoldings.

\subsection{General properties}

From now on we restrict our attention to complexes~$K$ satisfying the following two conditions:
\begin{enumerate}
\item
  $K$ is a triangulation of a PL-manifold;
\item
  $\odd{K}$ is a locally flat codimension-$2$-submanifold of $K$.
\end{enumerate}

It follows that both the unfolding and each component of the partial unfolding of~$K$ are PL-manifolds.

Making use of the homomorphism~$\frh_K$ defined in~(\ref{eq:frh}) we can reformulate the above question as follows.

\begin{prob}\label{ClassifyUnfoldings}
  Suppose that $N$ is a $d$-dimensional PL-manifold, $L$ is a locally flat codimension-$2$-submanifold, and
  $\frh:\pi_1(N\setminus L)\to S_{d+1}$ is a homomorphism.
  
  What are necessary and sufficient conditions such that there is a triangulation~$K$ of~$N$ with $\odd{K}=L$ and
  $\frh_K=\frh$?
\end{prob}

The following gives necessary conditions.

\begin{prp} \label{prp:PropertiesOfUnfolding}
  If such a triangulation~$K$ exists then:
  \begin{enumerate}
  \item the submanifold $L\subset N$ is a boundary~$\mod 2$;
  \item the homomorphism~$\frh$ takes each standard generator of $\pi_1(N\setminus L)$ to a transposition.
  \end{enumerate}
\end{prp}

\begin{proof}
  The first condition means that for some (and thus for any) triangulation of the pair $(N,L)$ there is a pure
  codimension-$1$-subcomplex~$Q$ of~$N$ such that~$\mod 2$ the formal sum of the boundaries of its facets equals the sum
  of the facets of~$L$.  In other words, the fundamental cycle of the submanifold~$L$ equals~$0$ in the homology
  group~$H_{d-2}(N,\ZZ_2)$.  If $K$ is a triangulation of a PL-manifold, then its odd subcomplex~$\odd{K}$ is the
  boundary~$\mod 2$ of the codimension-$1$-skeleton of~$K$.
  
  The term \emph{standard generator} is borrowed from the Wirtinger representation of the group of a knot.  A standard
  generator is the homotopy class of a loop~$\gamma$ in~$N\setminus L$, where $\gamma=\beta\lambda\beta^{-1}$
  and~$\lambda$ runs along the boundary of a small disk whose center lies on~$L$, and which is transversal to~$L$.  Such
  a disk exists due to our assumption on local flatness.  The projectivity along the discretization of such a loop (or
  of a loop close to it) is just a projectivity around a codimension-$2$-face.  Hence the proposition.
  
  Observe that if $N$ is not simply connected then the standard generators do not generate the fundamental group of the
  complement of~$L$.
\end{proof}

In the following section we will prove that these conditions are sufficient if $d=3$.  In particular, in this case there
are no restrictions on non-standard generators.

The following construction proves that a certain class of $2$-fold branched coverings can be realized as the complete
unfolding of a triangulation.

\begin{prp}
  Let $N$ be a PL-manifold and let $L$ be a codimension-$2$-sub\-poly\-he\-dron which is a boundary~$\mod 2$.  Then
  there exists a triangulation~$K$ of~$N$ such that $\odd{K}=L$ and $\frh_K[\gamma]=t^{l(\gamma,L)}$, where $t\in
  S_{d+1}$ is a transposition and $l(\gamma,L)$ is the linking number~$\mod 2$ of $\gamma$ and~$L$.
\end{prp}

Here $\gamma$ is supposed to be in general position with respect to~$Q$, where $Q$ is a codimension-$1$-subpolyhedron
whose boundary~$\mod 2$ is~$L$.  Then the \emph{linking number~$\mod 2$} of $\gamma$ and~$L$ is defined as the parity of
the number of intersections of~$\gamma$ and~$Q$.

\begin{proof}
  Take a triangulation~$(K',J')$ of the pair~$(N,L)$ such that $J'$ is an induced subcomplex of~$K'$.  Let $Q$ be a
  subcomplex of~$K'$ whose boundary~$\mod 2$ is~$J'$.  For each facet~$\tau$ of~$Q$ choose $\sigma$ to be one of the two
  facets of~$K'$ which contain~$\tau$.  By $e_\tau$ denote the edge in the barycentric subdivision~$b(K')$ of~$K'$ which
  connects the barycenters of $\sigma$ and~$\tau$.  We have the equality
  $$
  b(J')=\sum\lk_{b(K')}e_\tau\ (\mod 2),
  $$
  where the sum ranges over all facets~$\tau$ of~$Q$.  Now consider the complex~$K$ obtained from~$b(K')$ by
  subdivision of all edges~$e_\tau$.  Since $J'$ is an induced subcomplex of~$K'$, these subdivisions are independent of
  each other.  It is easy to see that $\odd{K}=b(J')$.
  
  Further, the vertices of~$b(K')$ can be colored with $d+1$~colors.  Let $0$ and~$1$ be the colors of the barycenters
  of facets and ridges of~$K'$, respectively.  Then the projectivity along a path~$\gamma$ is the $q$-th power of
  transposition~$(0\ 1)$, where $q$ counts how often $\gamma$ pierces~$Q$.
\end{proof}

\subsection{Unfoldings of $3$-dimensional manifolds}
\label{sec:3Unfolding}

In this section we prove the Characterization Theorem formulated in Section~\ref{sec:formulation}. 

It is known that each $1$-dimensional PL-submanifold (not necessarily connected, that is, a knot or a link) is
necessarily locally flat.

In fact, the monodromy homomorphism of a branched covering with the above properties is a homomorphism
$\frh:\pi_1(N\setminus L)\to S_4$ satisfying the conditions from Proposition~\ref{prp:PropertiesOfUnfolding}.  Therefore
it suffices to prove that there exists a triangulation~$K$ of~$N$ such that $\odd{K}=L$ and $\frh_K=\frh$.

The proof is organized as follows.  In the first step we construct a suitable triangulation of a regular neighborhood
of~$L$ in~$N$.  Then we extend the triangulation to the rest of the manifold, using handlebody decomposition.

\subsubsection*{Triangulation of a regular neighborhood}

First suppose that $L$ is connected.  Then $L\approx\Sph^1$.  Let $R$ be some regular neighborhood of~$L$ in~$N$, that
is, $R$ is a $3$-dimensional submanifold of~$N$ which geometrically collapses to~$L$; see Glaser~\cite[Vol.~I,
III.B]{Gla}.  In particular, $L$ is a strong deformation retract of~$R$. As shown in Moise~\cite[Chap.~24,
Theorem~11]{Moise}, the manifold~$R$ is PL-homeomorphic either to the solid torus~$T$, or to the solid Klein bottle~$F$.
Thus $R\setminus L$ is homotopy equivalent to the corresponding surface.

\begin{figure}[htbp] 
  \begin{center}
    \epsfig{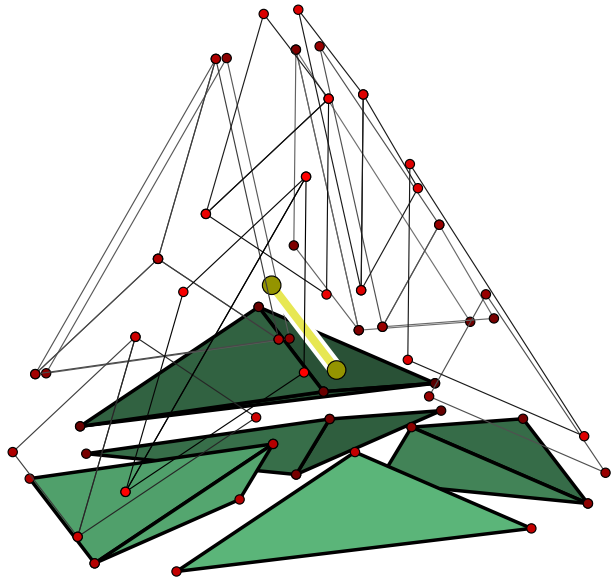}
    \caption{Explosion of the triangulated prism~$C_k$.  This complex has a
      rotational symmetry of order~$3$.  Only one fundamental domain is displayed as a set of five solid tetrahedra.}
    \label{fig:block}
  \end{center}
\end{figure}

Let us first consider the orientable case, that is, $R\setminus L\sim T$.  We start by analyzing the possible structure
of the homomorphism $\frh\circ i_*:\pi_1(R\setminus L)\to S_4$, where $i_*$ is the homomorphism induced by the inclusion
$i:R\setminus L\to N\setminus L$.  Let $a$ be the element of $\pi_1(R\setminus L)$ defined by a \emph{meridional loop}
($a$ is defined up to taking the inverse).  Pick an element $b\in \pi_1(R\setminus L)$ so that $a$ and $b$ together
generate $\pi_1(R\setminus L)$.  Since $i_*(a)$ is a standard generator in~$\pi_1(N\setminus L)$, that is, `a loop
around~$L$,' the permutation~$\frh\circ i_*(a)$ is a transposition, say $(0\ 1)$.  Since $a$ and $b$ commute, there are
four possible values for $\frh\circ i_*(b)$: either $\id$, $(0\ 1)$, $(2\ 3)$, or $(0\ 1)(2\ 3)$.  The second and the
third possibilities can be reduced to the first and the fourth one, respectively, by replacing $b$ with $ab$.  In order
to construct a suitable triangulation, subdivide $R$ into $n$~cylinders $C_1,C_2,\dots,C_n$, where $n$ is even if
$\frh\circ i_*(b)=\id$, and $n$ is odd if $\frh\circ i_*(b)=(0\ 1)(2\ 3)$.  We represent the cylinder~$C_k$ as a
triangular prism with triangular faces $x_{k-1}\, y_{k-1}\, z_{k-1}$ and $x_k\, y_k\, z_k$, respectively; indices are
taken modulo~$n$.  We assume that the closed PL-path $(x_0,x_1,\dots,x_{n-1},x_n=x_0)$ represents the element $b$~in
$\pi_1(R\setminus L)$.  Further, we assume that the intersection $C_k\cap L$ is an interval~$[v_{k-1},v_k]$, where the
point $v_k$ lies inside the triangle $x_k\, y_k\, z_k$.  Then we subdivide the prism $x_{k-1}\, y_{k-1}\, z_{k-1}\,
x_k\, y_k\, z_k$ into three triangular prisms with the common edge $v_{k-1}\,v_k$.  The resulting prism $x_{k-1}\,
y_{k-1}\, v_{k-1}\, x_k\, y_k\, v_k$ is triangulated into five tetrahedra
\begin{eqnarray*}
&\{v_{k-1},v_k,x_{k-1},y_{k-1}\},\\
&\{v_k,r_k,x_{k-1},y_{k-1}\}, \{v_k,r_k,y_{k-1},y_k\}, \{v_k,r_k,y_k,x_k\}, \{v_k,r_k,x_k,x_{k-1}\},
\end{eqnarray*}
where $r_k$ is an additional vertex inside the face $x_{k-1}\,y_{k-1}\,y_k\,x_k$.  In the same way we triangulate the
other two prisms $$y_{k-1}\, z_{k-1}\, v_{k-1}\, y_k\, z_k\, v_k\qquad\text{and}\qquad x_{k-1}\, z_{k-1}\, v_{k-1}\,
x_k\, z_k\, v_k$$ with additional vertices $s_k$ inside the face $y_{k-1}\,z_{k-1}\,y_k\,z_k$ and $t_k$ inside the face
$x_{k-1}\,z_{k-1}\,x_k\,z_k$, respectively.  For an illustration see Figure~\ref{fig:block}.

This way we obtain a triangulation of~$R$ which has the closed path $(v_0,v_1,\dots,v_{n-1},v_n=v_0)$ as its odd
subcomplex.  Now an arbitrary projectivity from $\{x_k,y_k,v_k,v_{k+1}\}$ to $\{x_l,y_l,v_l,v_{l+1}\}$ maps the pair
$(x_k,y_k)$ to either $(x_l,y_l)$ or $(y_l,x_l)$.  Moreover, it maps the pair $(v_k,v_{k+1})$ to $(v_l,v_{l+1})$ if and
only if $l-k$ is even; that is, $(v_k,v_{k+1})$ is mapped to $(v_{l+1},v_l)$ if $l-k$ is odd.

We choose $\sigma_0=\{x_0,y_0,v_0,v_1\}$ as our base facet and we identify $x_0\leftrightarrow 0$, $y_0\leftrightarrow
1$, $v_0\leftrightarrow 2$, and $v_1\leftrightarrow 3$.  Then the computation above implies that the projectivity along
any path whose homotopy class in $\pi_1(R\setminus L)$ equals~$b$ is either the identity or the double transposition
$(0\ 1)(2\ 3)$ provided that $n$ is even, and it equals either $(0\ 1)$ or~$(2\ 3)$ if $n$ is odd.

In the case that $R\setminus L\sim F$ we proceed exactly the same way.  Here the generators $a$ and $b$ satisfy the
relation $a^{-1}b=ba$.  Since the image of $a$ in $S_4$ is a transposition, again the images of $a$ and $b$ commute.
Compared to the case above the triangulation of $R$ differs only in that we identify $x_0$ with $y_n$ and $y_0$ with
$x_n$.  It is readily seen that again we can realize $\frh\circ i_*(b)$ as a projectivity by making $n$ either even
or odd.

Finally, if $L$ is not connected, then the regular neighborhoods of its components can be assumed disjoint. Their
triangulations are then constructed independently.

\subsubsection*{Extension of the triangulation to~$N\setminus R$}

Consider a relative handlebody decomposition of the pair $(N,R)$:
$$
R=N_{-1}\subset N_0\subset N_1\subset N_2\subset N_3=N,
$$
where $N_k$ is obtained from $N_{k-1}$ by attaching of a finite number of $k$-handles; see
Glaser~\cite[Vol.~II, p.~49]{Gla}.  This means that for each $k\in\{0,1,2,3\}$ we have a finite
collection of PL-embeddings $\cF_k=\{f_{k,i}:\partial \DD_i^k\times \DD_i^{3-k}\to\partial
N_{k-1}\}$ (where $\partial \DD^0=\emptyset$) with pairwise disjoint images.  We have
$$
N_k = N_{k-1}{\cup}_{\cF_k}\bigsqcup_i(\DD_i^k\times \DD_i^{3-k}),
$$

We will successively construct triangulations of the manifolds~$N_k$ with the following properties:

\begin{enumerate}
\item \label{cond:even} The odd subcomplex of $N_k$ is $L$.
\item \label{cond:hom} The homomorphism $\pi_1(N_k\setminus L)\to S_4$ defined by the triangulation
  as in~(\ref{eq:frh}) coincides with the composition of~$\frh$ with the homomorphism
  $\pi_1(N_k\setminus L)\to\pi_1(N\setminus L)$ induced by the inclusion $N_k\subseteq N$.
\end{enumerate}

Attaching $0$-handles is trivial: They form a collection of $3$-disks disjoint from~$R$, which can
be triangulated in such a way that their odd subcomplexes are empty.  We call a triangulation with
this property \emph{even}.  For example, take the barycentric subdivision of any triangulation.

In order to attach the $1$-handles suppose that the images of the attaching maps $\{f_{1,i}:\partial
\DD_i^1\times \DD_i^2\to\partial N_0\}$ are subcomplexes of $\partial N_0$.  This can be achieved by
applying iterated anti-prismatic subdivisions to~$N_0$, see
Proposition~\ref{prp:anti-prismatic-fine}.  Thus for each handle $\DD_i^1\times\DD_i^2$ the part
$\partial\DD_i^1\times\DD_i^2$ of its boundary is already triangulated and our aim is to extend this
triangulation in an appropriate way.  Fix an arbitrary collection of $3$-simplices
$\sigma_{i0},\sigma_{i1}$ in $N_0$, where $\sigma_{it}$ is adjacent to the component
$\{t\}\times\DD_i^2$ of $\partial\DD_i^1\times\DD_i^2=\{0,1\}\times\DD_i^2$.  Then for any even
triangulation of the handle~$\DD_i^1\times\DD_i^2$ the `projectivity along the handle'
$V(\sigma_{i0})\to V(\sigma_{i1})$ does not depend on the facet path chosen.

We show that, for any bijection $\phi_i:V(\sigma_{i0})\to V(\sigma_{i1})$, there is a triangulation
of the $i$-th handle~$\DD_i^1\times\DD_i^2$ such that the corresponding projectivity along the
handle coincides with~$\phi_i$.  For this purpose, fix a coloring of $\sigma_{i0}$.  It induces a
proper coloring of the simplicial neighborhood $$\bigcup\SetOf{\sigma}{\sigma\cap
  f_{1,i}(\{0\}\times\DD_i^2)\ne\emptyset}$$ of the subcomplex $f_{1,i}(\{0\}\times\DD_i^2)$ in~$N_0$.
Moreover, via~$\phi_i$, it also induces a proper coloring of the simplicial neighborhood of
$f_{1,i}(\{1\}\times\DD_i^2)$. In particular, we get a $4$-coloring of the $2$-dimensional complex
$\partial\DD_i^1\times\DD_i^2$.  Clearly, the existence of the desired triangulation of the handle
$\DD_i^1\times\DD_i^2$ follows from the following lemma, which will be proved later.

\begin{lem} \label{lem:1Handles}
  Each $4$-colored triangulation of $\partial\DD^1\times\DD^2$ extends to a $4$-colored
  triangulation of~$\DD^1\times\DD^2$.
\end{lem} 

Such triangulations of the $1$-handles extend the triangulation of~$N_0$ to a triangulation of~$N_1$
with the property~(\ref{cond:even}) above: There are no new odd edges since the simplicial
neighborhood of each handle is $4$-colorable.

Thus we obtain a triangulation of the manifold~$N_1$ with any prescribed projectivities along the
handles.  But, it is not hard to see that there is a collection $\{\phi_i:V(\sigma_{i0}\to
V(\sigma_{i1})\}$ of bijections such that the following holds: Any triangulation of~$N_1$ with
property~(\ref{cond:even}) which extends a triangulation of~$N_0$ and which realizes the
bijections~$\{\phi_i\}$ as projectivities along the handles satisfies the property~(\ref{cond:hom})
above.

We proceed to the $2$-handles.  Reasoning as before we can assume that
$B_i=f_{2,i}(\partial\DD_i^2\times\DD_i^1)$ is a subcomplex of~$N_1$.  The group of projectivities
of the simplicial neighborhood of~$B_i$ in~$N_1$ is trivial.  This follows from the fact that $B_i$
is contractible in~$N\setminus L$ and the property~(\ref{cond:hom}) of~$N_1$.  We conclude that the
simplicial neighborhood of~$B_i$ is $4$-colorable.  Hence in order to obtain a suitable
triangulation of~$N_2$ it is sufficient to prove the following lemma.

\begin{lem} \label{lem:2Handles}
  Each $4$-colored triangulation of~$\partial\DD^2\times\DD^1$ extends to a $4$-colored
  triangulation of~$\DD^2\times\DD^1$.
\end{lem}

Finally, the required triangulation of the $3$-handles is provided by the following result of
Goodman and Onishi~\cite{GO}.

\begin{lem}\label{lem:3Handles}
  Each $4$-colored triangulation of~$\partial D^3$ extends to a $4$-colored triangulation of $D^3$.
\end{lem}

This statement follows from Theorem~2.3 in~\cite{GO} (see the paragraph next to the theorem); it was
independently announced by Edwards~\cite{Edw}.

It remains to prove the Lemmas~\ref{lem:1Handles} and~\ref{lem:2Handles}.  This will be achieved by
a reduction to the Lemma~\ref{lem:3Handles}.

Suppose we have a $4$-colored triangulation of $\partial\DD^1 \times \DD^2$.  Take a
PL-homeo\-mor\-phism from $(\{0\} \times \DD^2)\setminus\relint\sigma$, where $\sigma$ is a facet in
the interior of~$\{0\} \times \DD^2$, onto the tube~$[0,\frac{1}{3}] \times \partial \DD^2$ such
that its restriction to~$\{0\}\times\DD^2$ is the identity.  This gives us an extension of the
$4$-colored triangulation to~$[0,\frac{1}{3}] \times \partial \DD^2$.  Similarly, we extend the
triangulation to the tube~$[\frac{2}{3},1] \times \partial \DD^2$ using the triangulation
of~$\{1\}\times\DD^2$.  Each of the circles $\{\frac{1}{3}\}\times\partial\DD^2$ and
$\{\frac{2}{3}\}\times\partial\DD^2$ is triangulated with exactly three vertices.  A simple
consideration shows that the $4$-colored triangulation can be extended to the
tube~$[\frac{1}{3},\frac{2}{3}] \times \partial \DD^2$.  Again we are in the position of
Lemma~\ref{lem:3Handles}.

As for Lemma~\ref{lem:2Handles}, the arguments in the first part of~\cite{GO} show that any
$4$-colored triangulation of the circle extends to a $4$-colored triangulation of the disk.  Hence
the given $4$-colored triangulation of $\DD^1\times\partial\DD^2$ extends to
$\partial(\DD^1\times\DD^2)$.  Again we are in the position of Lemma~\ref{lem:3Handles}.

This completes the proof of Theorem~\ref{thm:3dim}.

\subsection{Unfoldings of Surfaces}

Following the same line of reasoning one can prove the $2$-dimensional analogue of Theorem~\ref{thm:3dim}.

\begin{thm}\label{thm:2dim}
  Let $N$ be a closed surface, and let $f:M\to N$ be a branched covering with the following
  properties:
  \begin{enumerate}
  \item the number of sheets is less than or equal to~$3$;
  \item the number of branch points is even;
  \item the index of branching at any point in the pre-image of a branch point is either~$1$ or~$2$.
  \end{enumerate}
  Then there is a triangulation~$K$ of~$N$ such that $M$ is PL-homeomorphic to a component of the
  partial unfolding of~$K$ and $f$ is equivalent to the restriction of the partial unfolding map.
\end{thm}

We give a direct proof for the $2$-dimensional analogue to Theorem~\ref{thm:main-result}.

\begin{thm}
  For each closed orientable surface~$M_g$ of genus~$g$ there is a triangulation~$P_g$ of the
  sphere~$\Sph^2$ such that one of the components of its partial unfolding is homeomorphic to~$M$.
  Moreover, this component is a $2$-fold branched covering, and hence is isomorphic to the
  complete unfolding of~$P_g$.
\end{thm}

\begin{proof}
  The Riemann-Hurwitz formula expresses the Euler characteristic of a branched covering space~$M$
  over the $2$-sphere in terms of the number of sheets and the branching indices.  In particular, if
  $f:M_g\to \Sph^2$ is a $2$-fold branched covering with $2n$ branch points (thus the pre-image of
  each branch point consists of a single point of index~$2$) we have
  $$
  \chi(M)=2\chi(\Sph^2)-2n=4-2n.
  $$
  Hence it suffices to construct a triangulation $P_g$ of~$\Sph^2$ whose partial unfolding is a
  $2$-fold branched covering with exactly $2(g+1)$ branch points.  Actually, the existence of such a
  triangulation follows from Theorem~\ref{thm:2dim}. But, we present here an explicit construction.

  For $P_0$ one can take the suspension over the boundary of the triangle.
  
  Let $g>0$.  Then there exists a triangulation~$Q_g$ of~$\Sph^2$ with $2(g+1)$ facets.  Starting
  from the boundary of a simplex the triangulation~$Q_g$ can be constructed inductively by stellar
  subdivision.  Then perform a stellar subdivision on each facet of~$Q_g$.  We denote the resulting
  triangulation by~$P_g$.  It has has $2(g+1)$ vertices of degree~$3$, all the other vertices being
  even.  We have $\Pi(P_g)=\ZZ_2$.
  
  Note that the unfolding of~$P_g$ coincides with the partial unfolding with respect to any even
  vertex.
\end{proof}

\appendix\section{Appendix}\label{sec:TechLemmas}

We postponed some of the more technical details until now.

\subsection{Anti-prismatic subdivision}\label{subsec:anti-prismatic}

The group of projectivities and the unfoldings of a simplicial complex are invariants of its
combinatorial structure.  Different triangulations of the same topological space usually yield
different groups of projectivities.  In particular, the barycentric subdivision~$b(K)$ of a
simplicial complex~$K$ always has a trivial group of projectivities.  In view of
Proposition~\ref{prp:Identity} this implies that $\unfolding{b(K)}$ is isomorphic to~$b(K)$ for
locally strongly connected~$K$.

The barycentric subdivision plays a crucial role in many technical aspects of PL-topology.  As its
key feature the iterated barycentric subdivision becomes arbitrarily fine.  Here we need a fine
subdivision which respects the group of projectivities.  As already pointed out, the barycentric
subdivision is out of question.  Therefore, as an alternative, we suggest the \emph{anti-prismatic
  subdivision}.

Let $c_n$ denote the simplicial complex arising from the Schlegel diagram of (the boundary of) the
$(n+1)$-dimensional cross polytope.  More precisely, $c_n$ has $2n+2$ vertices
$v_0^+,\dots,v_n^+,v_0^-,\dots,v_n^-$ and $2^{n+1}-1$~facets which correspond to the subsets~$S$ of
the vertices with $n$~elements such that $S$ contains exactly one vertex from each antipodal pair
$\{v_i^+,v_i^-\}$; the set $\{v_0^+,\dots,v_n^+\}$ is excluded.  Then we can consider $c_n$ as a
subdivision of the geometric simplex with vertices $v_0^+,\dots,v_n^+$.  For an introduction to
Schlegel diagrams of convex polytopes see Ziegler~\cite[Chapter~5]{GMZ}.

Now suppose that $\tau=\{v_0,\dots,v_n\}$ is some face of the complex~$K$.  The operation of
\emph{crossing} of the face $\tau$ in~$K$ replaces the star of~$\tau$ by the join of~$c_n$ with the
link of $\tau$:
$$
c(K,\tau)=(K\setminus\st\tau)\cup(c_n*\lk\tau).
$$
Since $c_n\approx\tau$, the result of crossing is PL-homeomorphic to the initial complex~$K$.  In
order to obtain the \emph{anti-prismatic subdivision}~$a(K)$ of $K$ the crossings of the faces of
$K$ must be performed in all dimensions with decreasing order: We start with the facets and then go
down to edges.  Observe that the crossing of a vertex is trivial.

The anti-prismatic subdivision is analogous to the barycentric subdivision in the sense that the
crossing in the former plays the same role as the \emph{starring} in the latter.  While crossing
means to substitute a $k$-simplex by the diagram of the $(k+1)$-dimensional cross-polytope, starring
substitutes a $k$-simplex by the diagram of the $(k+1)$-dimensional simplex.

\begin{prp}\label{prp:anti-prismatic-fine}
  Let $K$ be any geometric simplicial complex.  For each $\varepsilon>0$ there is a natural
  number~$n$ such that each simplex of an $n$-times iterated anti-prismatic subdivision of~$K$ has
  diameter less than~$\varepsilon$.
\end{prp}

This can be proved in the same way as the corresponding result on the barycentric subdivision.

The anti-prismatic subdivision can also be defined for pseudo-simplicial complexes.  We have the
following central property of the anti-prismatic subdivision.

\begin{prp}\label{prp:anti-prismatic-simplicial}
  The anti-prismatic subdivision of a pseudo-simplicial complex is a simplicial complex.
\end{prp}

\begin{proof}
  We give an alternative description of~$a(K)$ as an abstract simplicial complex.
  
  As the set of vertices take all pairs $(\tau,w)$, where $\tau$ is a non-empty face of~$K$, and $w$
  is a vertex of~$\tau$.  An original vertex~$w$ of~$K$ is naturally identified with the
  vertex~$(w,w)$ of~$a(K)$.  The set $\{(\tau_0,w_0),\ldots,(\tau_k,w_k)\}$ is a face of the
  anti-prismatic subdivision if and only if
  \begin{enumerate}
  \item $\tau_0\subseteq\tau_1\subseteq\cdots\subseteq\tau_k$ is a flag in~$K$ (with repetitions
    allowed), and
  \item if $\tau_i\subsetneq\tau_j$ then $w_j\notin\tau_i$.
  \end{enumerate}
  Provided that all the pairs $(\tau_0,w_0),\ldots,(\tau_k,w_k)$ are distinct, or, equivalently, the
  face~$\{(\tau_0,w_0),\ldots,(\tau_k,w_k)\}$ is of dimension~$k$, it follows that all the
  vertices~$w_0,\ldots,w_k$ are distinct.

  \begin{figure}[htbp]
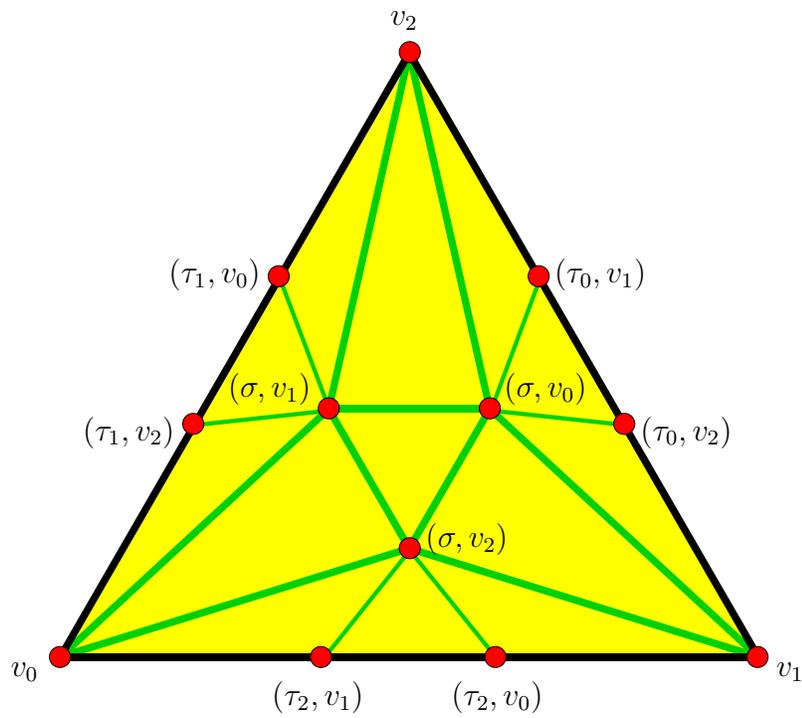

    \begin{center}
      \input anti-prismatic.pstex_t 
      \caption{Anti-prismatic subdivision of the triangle $\sigma$ with vertices $v_0,v_1,v_2$ and
        edges $\tau_0,\tau_1,\tau_2$.}
      \label{fig:anti-prismatic}
    \end{center}
  \end{figure}
  
  Use an induction to see that this yields the same as the construction by iterated crossing.  The
  vertices $(\sigma,v_0),\ldots,(\sigma,v_d)$ emerge as the result of the crossing of the
  $n$-dimensional face~$\sigma$ of~$K$ with $V(\sigma)=\{v_0,\ldots,v_d\}$; here $v_i$ is opposite
  to~$(\sigma,v_i)$, see Figure~\ref{fig:anti-prismatic}.
\end{proof}

In the following we use the description of~$a(K)$ which was given in the proof of
Proposition~\ref{prp:anti-prismatic-simplicial}.

The map $f\colon(\tau,w)\mapsto w$ is a non-degenerate simplicial map from $a(K)$ onto~$K$.  We call
$f$ the \emph{crumpling map} of the anti-prismatic subdivision.

As an immediate consequence the $k$-colorability of the $1$-skeleton of~$K$ implies the
$k$-colorability of the $1$-skeleton of~$a(K)$.  Besides, $K$ is balanced if and only if
$a(K)$ is balanced.  However, a stronger property holds.

By Proposition~\ref{prp:non-degenerate} we obtain a monomorphism~$f_*$ between the groups of
projectivities.

\begin{prp}\label{prp:isomorphic-proj-grp}
  The induced map $$f_*:\Pi(a(K))\to\Pi(K)$$ is an ismorphism of permutation groups.
\end{prp}

\begin{figure}[htbp]
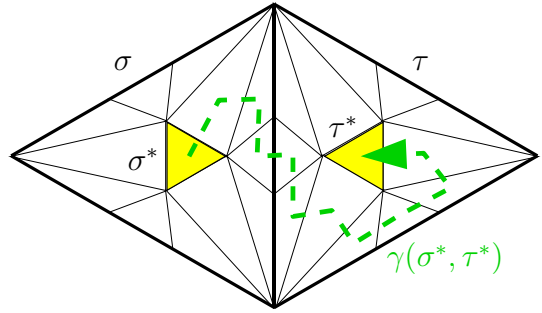

  \begin{center}
    \input anti-prismatic-2.pstex_t 
    \caption{Perspectivity $\langle\sigma,\tau\rangle$ in~$K$ lifted to the 
      projectivity $\langle\gamma(\sigma^*,\tau^*)\rangle$ in~$a(K)$.}
    \label{fig:anti-prismatic-2}
  \end{center}
\end{figure}

\begin{proof}
  We have to verify that $f_*$ is surjective.  As a stronger property we actually show that each
  perspectivity in~$K$ can be `lifted' to a projectivity in~$a(K)$.  For each facet~$\sigma$
  of~$K$ let $$\sigma^*=\SetOf{(\sigma,v)}{v\in V(\sigma)}$$ be the \emph{corresponding} facet
  of~$a(K)$.  Clearly, $f(\sigma^*)=\sigma$.  In order to lift the
  perspectivity~$\langle\sigma,\tau\rangle$ choose an arbitrary facet path~$\gamma(\sigma^*,\tau^*)$
  from $\sigma^*$ to~$\tau^*$ in the subcomplex $a(\sigma\cup\tau)$ of~$a(K)$.  Since the complex
  $a(\sigma\cup\tau)$ is balanced, the resulting projectivity
  $\langle\gamma(\sigma^*,\tau^*)\rangle:V(\sigma^*)\to V(\tau^*)$ does not depend on the choice of
  the facet path~$\gamma(\sigma^*,\tau^*)$; see Figure~\ref{fig:anti-prismatic-2}.  For any facet
  path $(\sigma_0,\sigma_1,\ldots,\sigma_n=\sigma_0)$ we obtain
  $$f_*(\langle\gamma(\sigma_0^*,\sigma_1^*)\rangle\cdots\langle\gamma(\sigma_{n-1}^*,\sigma_0^*)\rangle)
  =
  \langle\sigma_0,\sigma_1\rangle\cdots\langle\sigma_{n-1},\sigma_0\rangle.$$
\end{proof}

\begin{prp}\label{prp:trillian}
  The unfolding $\unfolding{a(K)}$ is canonically isomorphic to $a(\unfolding{K})$ as a pseudo-simplicial
  complex.  In particular, $\unfolding{a(K)}$ is a simplicial complex.
\end{prp}

\begin{proof}
  We start by scrutinizing the construction of the unfolding from Section~\ref{subsec:unfolding}.
  It is essential that the unfolding of a complex~$K$ is obtained from gluing the geometric
  simplices $(\sigma,g)$, where $\sigma$ is a facet of~$K$ and $g$ is a projectivity, in an
  \emph{arbitrary} order.
  
  The anti-prismatic subdivision of a simplex is locally strongly connected and balanced.  From
  Proposition~\ref{prp:Identity} we infer that its unfolding is isomorphic to itself.  Therefore,
  for each facet $\sigma$ of~$K$, we can first glue the facets in~$a(\sigma)$.  In the second step
  we glue these subdivided facets~$a(\sigma)$ to obtain $\unfolding{a(K)}$.  This is equivalent to
  the construction of~$\unfolding{K}$, because, for each vertex~$(\tau,v)$ of~$a(K)$ with
  $\codim\tau=1$ the star~$\st_{a(K)}(\tau,v)$ is balanced and~$\Pi(a(K))\cong\Pi(K)$.
\end{proof}

\begin{cor}\label{cor:diagram1}
  The following diagram of pseudo-simplicial complexes and simplicial maps is commutative.
  $$
  \begin{array}{ccc}
    \unfolding{a(K)} & \rightarrow & \unfolding{K} \\[0.5ex]
    \downarrow       &             & \downarrow\\[0.5ex]
    a(K)             & \rightarrow & K
  \end{array}
  $$
  The vertical arrows are complete unfolding maps, while the horizontal arrows are the crumpling
  maps.
\end{cor}

\begin{cor}\label{cor:diagram2}
  The following diagram of topological spaces and continuous maps is commutative.
  $$
  \begin{array}{ccc}
    \unfolding{a(K)} & \rightarrow & \unfolding{K} \\[0.5ex]
    \downarrow       &             & \downarrow\\[0.5ex]
    a(K)             & \rightarrow & K
  \end{array}
  $$
  The vertical arrows are complete unfolding maps, while the horizontal arrows are PL-homeomorphisms
  induced by subdivision.
\end{cor}

For the partial unfolding the situation is completely analogous.  One can prove the following.

\begin{prp}
  The unfolding $\partunfold{a(K)}$ is canonically isomorphic to $a(\partunfold{K})$ as a pseudo-simplicial
  complex.  In particular, $\partunfold{a(K)}$ is a simplicial complex.
\end{prp}

We obtain commutative diagrams which are similar to the ones in Corollary~\ref{cor:diagram1} and
Corollary~\ref{cor:diagram2}.

\subsection{Homotopy properties of nice complexes}\label{subsec:homotopy}

Here we prove that locally strong connectivity and locally strong simple connectivity (as defined in
Section~\ref{sec:Odd}) provide good homotopy properties of the dual skeleta of the simplicial
complex.  Namely, they allow to approximate paths and homotopies by paths and homotopies in the dual
$1$-skeleton and the dual $2$-skeleton, respectively.

Note that there is a somewhat similar situation for the homology properties: A \emph{triangulated
  homology manifold} is a simplicial complex such that the link of each face is a homology sphere of
the appropriate dimension.  This property implies a Poincar\'e duality theorem, see
Munkres~\cite[\S65]{Mun}.  Thus the (local) homology properties of the links provide a good (global)
homology structure for the whole complex.

Again $K$ is a pure and locally finite simplicial complex.

Recall that $b(K)$ denotes the barycentric subdivision of a simplicial complex~$K$, and the
simplices of~$b(K)$ have the form $\{\hat\sigma_0,\ldots,\hat\sigma_n\}$, where
$\sigma_0\subset\dots\subset\sigma_n$ are faces of~$K$; the point~$\hat\sigma_i$ is the barycenter
of the face~$\sigma_i$.  If $\sigma$ is a face of~$K$ then the \emph{block dual to}~$\sigma$ is the
geometric subcomplex of~$b(K)$ defined by
\begin{equation} \label{eqn:BlockDefinition}
D(\sigma)=\bigcup_{\sigma_0=\sigma}\conv\{\hat\sigma_0,\ldots,\hat\sigma_n\}.
\end{equation}
The interior of the block is
$$
\int D(\sigma)
=\bigcup_{\sigma_0=\sigma}\relint(\conv\{\hat\sigma_0,\ldots,\hat\sigma_n\})
$$
and its boundary is
$$
\partial D(\sigma)=\bigcup_{\sigma_0\subsetneq\sigma}\conv\{\hat\sigma_0,\ldots,\hat\sigma_n\}.
$$

It is readily seen that $D(\sigma)=\partial D(\sigma) \sqcup\int D(\sigma)$ and $K
=\bigsqcup_{\sigma\in K}\int D(\sigma)$. Note also that
\begin{equation} \label{eqn:BlockAsCone}
D(\sigma)=\hat\sigma*\partial D(\sigma)\approx\cone\partial D(\sigma).
\end{equation}
The block $D(\sigma)$ is sometimes called the \emph{barycentric star} of~$\sigma$ since
$D(\sigma)=\st_{b(K)}\hat\sigma$.  By the same token the boundary $\partial D(\sigma)$ of the block
is called the \emph{barycentric link} of~$\sigma$.

By the \emph{dual block complex} $K^*$ we mean the geometric realization of the complex $b(K)$
together with the decomposition $K=\bigcup_{\sigma\in K}D(\sigma)$.  Due to $K$ being pure each
block is a pure subcomplex.  We have
$$
\dim D(\sigma)=\codim\sigma.
$$
Now the space
$$
K^*_{(n)}=\bigcup_{\codim\sigma=n}D(\sigma)
$$
is called the \emph{$n$-dimensional dual skeleton}.

\begin{prp} \label{prp:PathApprox}
  Let $K$ be a locally strongly connected simplicial complex and $f:[0,1]\to K^*$ be a continuous
  map with $f(0),f(1)\in K^*_{(0)}$. Then there exists a map $g:[0,1]\to K^*$ homotopic to $f$ with
  fixed endpoints such that $g([0,1])\subset K^*_{(1)}$.
\end{prp}

\begin{proof}
  Suppose that $f([0,1])\subset K^*_{(n)}$, where $n>1$.  We will show that there is a homotopy with
  fixed endpoints which deforms $f$ to a map whose image is contained in~$K^*_{(n-1)}$.  Thus by
  induction we obtain the desired map~$g$.
  
  Since $K$ is locally finite and $[0,1]$ is compact, the image $f([0,1])$ is contained in a finite
  number of $n$-dimensional blocks.  Let $D=D(\sigma)$ be one of them.  We construct a homotopy
  `sweeping' $f([0,1])$ out of~$\int D$.
  
  The set $f^{-1}(\int D)$ is a disjoint union of open intervals whose endpoints lie in $\partial
  D$.  The homotopy sequence of the pair
  $$
  \cdots\to\pi_1(D)\to \pi_1(D,\partial D)\to \pi_0(\partial D)\to\cdots
  $$
  is exact (in the wider sense of maps of pointed sets as in Bredon~\cite[VII.5]{Bred}).  Since
  the left and the right terms are trivial, the middle term also vanishes.  Hence we can deform the
  restrictions of the map~$f$ over all intervals so that the images will lie in~$\partial D$.
\end{proof}

\begin{prp} \label{prp:BlockApprox}
  Let $K$ be a nice simplicial complex and let $f:\DD^2\to K^*$ be a continuous map with
  $f(\Sph^1)\subset K^*_{(1)}$.  Then there exists a map $g:\DD^2\to K^*$ which is
  $\Sph^1$-homotopic to $f$ and such that $g(\DD^2)\subset K^*_{(2)}$.
\end{prp}

\begin{proof}
  The proof is essentially the same as that of the Whitehead's original proof~\cite{Whi} of the
  Cellular Approximation Theorem.
  
  By an induction, as above, it suffices to `sweep' the map $f$ out of a maximal block $D=D(\sigma)$
  with $\dim\sigma>2$.
  
  Two cases are to be considered.  If $\hat\sigma\notin f(\DD^2)$, then we use the fact that
  $\partial D$ is a strong deformation retract of $D\setminus\{\hat\sigma\}$, see
  (\ref{eqn:BlockAsCone}).  Otherwise, we have $\hat\sigma\in f(\DD^2)$.  Our aim is then to `free'
  the point $\hat\sigma$ from the image $f(\DD^2)$ thereby returning to the previous case.
  Triangulate the disk $\DD^2$ finely enough so that for any face~$\kappa$ of the triangulation
  \begin{equation} \label{eqn:FineTriang}
    \hat\sigma\in f(\kappa)\ \text{implies that}\ f(\kappa)\subseteq B.
  \end{equation}
  We proceed by an induction on the dimension of the simplices in the triangulation of~$\DD^2$ whose
  images cover the point~$\hat\sigma$, starting from dimension~$0$.  Suppose that $\kappa$ is a face
  of the triangulation such that $\hat\sigma\in f(\relint\kappa)$ and $f(\partial\kappa)\subset
  D\setminus\{\hat\sigma\}$.  From the exact sequence of the pair it follows that the relative
  homotopy groups $\pi_i(D,D\setminus\{\hat\sigma\})$ are trivial for $i=0,1,2$.  Thus there is a
  homotopy $f_t|_\kappa$ of the map $f|_\kappa$ which is constant on $\partial\kappa$ and such that
  $f_1(\kappa)\subset D\setminus\{\hat\sigma\}$.  Put $f_t$ constant on the closure of the
  complement to~$\st\kappa$ and apply Borsuk's Homotopy Extension Theorem, see
  Bredon~\cite[VII.1.4]{Bred} in order to obtain a homotopy on the whole space~$K^*$.  In this manner
  we can consecutively free the point~$\hat\sigma$ from all the faces of the triangulation
  of~$\DD^2$.\\
  \nothing\hfill
\end{proof}
\break

\bibliographystyle{amsplain}
\bibliography{brcov}
\vfill

\noindent
Ivan Izmestiev\\
NWF-I Mathematik\\
Universit\"at Regensburg\\
93040 Regensburg, Germany\\
\texttt{ivan.izmestiev@mathematik.uni-regensburg.de}

\bigskip\noindent
Michael Joswig\\
Institut f\"ur Mathematik, MA 6-2\\
Technische Universit\"at Berlin\\
Stra\ss{}e des 17.~Juni 136\\
10623 Berlin, Germany\\
\texttt{joswig@math.tu-berlin.de}

\end{document}